\documentclass[12pt]{amsart}

\DeclareMathOperator{\GL}{GL} %
\DeclareMathOperator{\PGL}{PGL} %
\DeclareMathOperator{\SL}{SL} %
\DeclareMathOperator{\U}{U} %
\DeclareMathOperator{\SU}{SU} %
\DeclareMathOperator{\Par}{Par} %
\DeclareMathOperator{\Ind}{Ind} %
\DeclareMathOperator{\Irr}{Irr} %
\DeclareMathOperator{\fpr}{fpr} %
\DeclareMathOperator{\fix}{fix} %

\renewcommand\AA{{\mathbb A}}
\newcommand\QQ{{\mathbb Q}}
\newcommand\ZZ{{\mathbb Z}}
\newcommand\barQl{{\bar{\mathbb Q}_l}}
\newcommand\FF{\mathbb F}
\newcommand\NN{\mathbb N}

\newcommand\CB{{\mathcal B}}
\newcommand\CC{\mathcal C}
\newcommand\CF{{\mathcal F}}

\newcommand\CP{{\mathcal P}}
\newcommand\CQ{{\mathcal Q}}
\newcommand\CR{\mathcal R}
\newcommand\CS{{\mathcal S}}
\newcommand\CT{{\mathcal T}}
\newcommand\CW{{\mathcal W}}

\newcommand\uni{{\mathrm{uni}}}
\newcommand\uu{{\mathrm{u}}}

\newcommand\Flag{{\mathrm{Flag}}}

\newcommand\sse{\subseteq}
\renewcommand\bar{\overline}

\usepackage{amssymb}
\usepackage{fullpage}
\usepackage{color}

\numberwithin{equation}{section}

\theoremstyle{plain}
\newtheorem{lem}[equation]{Lemma}
\newtheorem{thm}[equation]{Theorem}
\newtheorem{cor}[equation]{Corollary}
\newtheorem{prop}[equation]{Proposition}

\theoremstyle{definition}

\newtheorem{exmp}[equation]{Example}
\newtheorem{exmps}[equation]{Examples}

\theoremstyle{remark}
\newtheorem{rem}[equation]{Remark}
\newtheorem{rems}[equation]{Remarks}

\title[Rational points and unipotent conjugacy]
{Rational points on generalized flag varieties and unipotent
conjugacy in finite groups of Lie type}
\author[S.~M.~Goodwin and G.~R\"ohrle]
{Simon M.~Goodwin and Gerhard R\"ohrle}

\address{School of Mathematics, University of Birmingham,
Birmingham, B15 2TT, United Kingdom}
\email{goodwin@maths.bham.ac.uk}
\urladdr{http://web.mat.bham.ac.uk/S.M.Goodwin}

\address
{Fakult\"at f\"ur Mathematik,
Ruhr-Universit\"at Bochum,
D-44780 Bochum, Germany}
\email{gerhard.roehrle@rub.de}
\urladdr{http://www.ruhr-uni-bochum.de/ffm/Lehrstuehle/Lehrstuhl-VI/rubroehrle.html}


\dedicatory{Dedicated to Professor J.~A.~Green on the occasion of
his 80th birthday}


\begin{document}

\begin{abstract}
Let $G$ be a connected reductive algebraic group defined over the
finite field $\FF_q$, where $q$ is a power of a good prime for $G$.
We write $F$ for the Frobenius morphism of $G$ corresponding to the
$\FF_q$-structure, so that $G^F$ is a finite group of Lie type. Let
$P$ be an $F$-stable parabolic subgroup of $G$ and $U$ the unipotent
radical of $P$.  In this paper, we prove that the number of
$U^F$-conjugacy classes in $G^F$ is given by a polynomial in $q$,
under the assumption that the centre of $G$ is connected. This
answers a question of J.~Alperin in \cite{alperin}.

In order to prove the result mentioned above, we consider, for
unipotent $u \in G^F$, the variety $\CP^0_u$ of $G$-conjugates of $P$
whose unipotent radical contains $u$. We prove that the number of
$\FF_q$-rational points of $\CP^0_u$ is given by a polynomial in $q$
with integer coefficients.
Moreover, in case $G$ is split over $\FF_q$ and $u$
is split (in the sense of \cite[\S5]{shoji}), the coefficients of
this polynomial are given by the Betti numbers of $\CP^0_u$.
We also prove the analogous results for the
variety $\CP_u$ consisting of conjugates of $P$ that contain $u$.
\end{abstract}

\maketitle

\section{Introduction}
\label{sect:intro}

Let $\GL_n(q)$ be the general linear group of nonsingular $n \times
n$ matrices over the finite field $\FF_q$ and let $\U_n(q)$ be the
subgroup of  $\GL_n(q)$ consisting of upper unitriangular matrices.
A longstanding conjecture states that the number of conjugacy
classes of $\U_n(q)$ for fixed $n$ as a function of $q$ is an
integral polynomial in $q$. This conjecture has been attributed to
G.~Higman (cf.\ \cite{higman}), and it has been verified for $n \le 13$, see
\cite{veralopezarregi}. There has also been interest in this
conjecture from G.R.~Robinson (see  \cite{robinson}) and J.~Thompson
(see \cite{thompson}).

In \cite{alperin} J.~Alperin showed that a related question is
easily answered, namely that the number of $\U_n(q)$-conjugacy
classes in all of $\GL_n(q)$, for fixed $n$ as a function of $q$ is a
polynomial in $q$ with integer coefficients. Generalizing this
result, in Theorem \ref{thm:main}, we give an affirmative answer to
a question raised by Alperin in \cite{alperin}.  In order to state
this generalization we need to introduce some notation.

Let $G$ be a connected reductive linear algebraic group
over $\FF_q$, where $q$ is a power of a good prime
for $G$.  We write $F$ for the
Frobenius morphism corresponding to the $\FF_q$-structure on $G$;
then the $F$-fixed point subgroup $G^F$ of $G$ is a finite group of
Lie type. Assume that the centre of $G$ is connected.
Let $P$ be an
$F$-stable parabolic subgroup of $G$ with unipotent radical $U$.

In Theorem \ref{thm:main} of this paper, we prove that the number of
$U^F$-conjugacy classes in $G^F$ is given by a polynomial in $q$ with
integer coefficients (in case $G$ has a simple component of type
$E_8$, we require two polynomials depending on the congruence of $q$
modulo $3$). This generalizes the aforementioned result of Alperin,
firstly by replacing $\GL_n(q)$ by an arbitrary finite reductive
group $G^F$ and secondly by allowing any parabolic subgroup of $G^F$
rather than just a Borel subgroup.

\smallskip

In order to prove Theorem \ref{thm:main}, we are led to consider the
$\FF_q$-rational points in certain subvarieties of the variety $\CP$
of $G$-conjugates of $P$. Given $u \in G$ unipotent, we write
$\CP_u$ for the subvariety of $\CP$ consisting of conjugates of $P$
that contain $u$ and $\CP^0_u$ is the subvariety of $\CP_u$
consisting of conjugates of $P$ whose unipotent radical contains
$u$. We prove that the number of $\FF_q$-rational points of the
variety $\CP_u$ (respectively $\CP_u^0$) is given by a polynomial in
$q$ for $u \in G^F$ unipotent.
For a precise formulation, see Theorem \ref{thmmain}. Assuming
that $G$ is split and that $u \in G^F$ is split (in the sense of
\cite[\S5]{shoji}) we show in Proposition \ref{p:betti} that the
coefficients of the polynomials describing the number of
$\FF_q$-rational points in $\CP_u$ and $\CP^0_u$ from Theorem
\ref{thmmain} are given by the Betti numbers of the underlying
varieties.
\medskip

We now give an outline of the proofs of our principal results
and of the structure of the paper.

In order to state our main results precisely, we require an axiomatic formalism
for connected reductive algebraic groups; this axiomatic setup is given
in \S \ref{s:setup}.  Further, to make sense of considering the varieties
$\CP_u$ for different values of $q$, we need a parameterization of the
unipotent $G^F$-conjugacy classes that is independent of $q$; the relevant
results from \cite{shoji} and \cite{kawanaka} are recalled in
\S \ref{s:classes}.

For our proof of Theorem \ref{thmmain},  we require results from
the representation theory of finite groups of Lie type to express
$|\CP_u^F|$ and $|(\CP_u^0)^F|$ in terms of Deligne--Lusztig
generalized characters. Consequently, $|\CP_u^F|$ and $|(\CP_u^0)^F|$
are given as linear combinations of Green functions.  Due to work of
T.~Shoji (\cite{shoji}) these Green functions are given by
polynomials in $q$; from this we can deduce Theorem \ref{thmmain}.
The relevant results on Deligne--Lusztig generalized characters and
Green functions are recalled in \S \ref{s:dlgenchars}.

In order to prove Proposition \ref{p:betti}, i.e.\ that, in case $G$
and $u$ are split, the coefficients of the polynomials $|\CP_u^F|$
and $|(\CP_u^0)^F|$ are given by the Betti numbers of the varieties
$\CP_u$ and $\CP_u^0$, we require a purity result for these varieties.
In \S \ref{s:purity} we deduce such a result from a well-known
purity result of T.~A.~Springer \cite{springer1} for the variety
$\CB_u$ of Borel subgroups of $G$ that contain $u$. This deduction
also requires results of W.~Borho and R.~MacPherson
(\cite{borhomac}) about the varieties $\CP_u$ and $\CP_u^0$.

By $k(U^F,G^F)$ we denote the number of $U^F$-conjugacy classes in
$G^F$. To prove Theorem \ref{thm:main} we adapt the arguments in
\cite{alperin} to express $k(U^F,G^F)$ in terms of the sizes of the
sets $(\CP_u^0)^F$ for $u \in U^F$. The independence of the
parameterization of the unipotent $G^F$-conjugacy classes means
that, in a sense, this expression does not depend on $q$. This
allows us to deduce Theorem \ref{thm:main} from Theorem
\ref{thmmain}.

Let $P$ and $Q$ be associated $F$-stable parabolic subgroups of $G$,
i.e.\ $P$ and $Q$ have Levi subgroups which are conjugate in $G$.
Let $\CP$ and $\CQ$ be the corresponding generalized flag varieties
of $P$ and $Q$, respectively. In Corollary \ref{cor6} we show that
for $u \in G^F$ unipotent, we have $|\CP_u^F| = |\CQ_u^F|$ and
$|(\CP_u^0)^F| = |(\CQ_u^0)^F|$. Further, let $U$ and $V$ be the
unipotent radicals of $P$ and $Q$, respectively. As a consequence of
the expressions that we get for the polynomials in Theorem
\ref{thm:main}, we show in Corollary \ref{cor2}, that the two
polynomials in $q$ given by $k(U^F,G^F)$ and $k(V^F,G^F)$ are equal.
This essentially follows from the fact that the Harish-Chandra
induction functor $R_L^G$ does not depend on the choice of an
$F$-stable parabolic subgroup of $G$ containing $L$ as a Levi
subgroup.

We note that the assumption that the centre of $G$ is connected is
only really required so that the parameterization of the unipotent
$G^F$-conjugacy classes is independent of $q$. This means that many
of the results in this paper are valid without this assumption given
an appropriate statement; for example, see Remark \ref{r:r16}.  In
Example \ref{ex:disc} we indicate in case $G$ has disconnected
centre how the number of $U^F$-conjugacy classes in $G^F$ is given
by a set of polynomials in $q$ depending on the congruence class of
$q$ with respect to some positive integer.

\smallskip

In an appendix to this paper, we give an elementary combinatorial
proof that $|(\CP_u^0)^F|$ is given by a polynomial in $q$ in the
case $G = \GL_n$.  Consequently, we obtain an elementary proof of
Theorem \ref{thm:main} in the case $G = \GL_n$.

\medskip

As general references on finite groups of Lie type, we refer the
reader to the books by Carter \cite{carter} and Digne--Michel
\cite{dignemichel}.

\section{Preliminaries}
\label{s:prelims}

For the statements of the main results of this paper we require an
axiomatic formalism for connected reductive algebraic groups; this is
introduced in \S \ref{s:setup}.  In other parts of the sequel we do
not require this axiomatic setup and we use the notation given
in the following subsection instead.
In each subsection, we state which
setup we are using.

\subsection{General notation for algebraic groups}
\label{s:notn}
In this section we introduce various pieces of notation that
we require in the sequel.

First we fix some general group theoretic notation.  For a group
$G$, a subgroup $H$ of $G$ and $x \in G$ we write: $Z(G)$ for the
centre of $G$; $N_G(H)$ for the normalizer of $H$ in $G$; $C_G(x)$
for the centralizer of $x$ in $G$; $G \cdot x$ for the conjugacy
class of $x$ in $G$; and ${}^xH = xHx^{-1}$. If $G$ is finite and
$p$ is a prime, then we write $|G|_p$ for the $p$-part of the order
of $G$, and $|G|_{p'}$ for the $p'$-part of the order of $G$.

Now we give some notation for connected reductive algebraic groups that we use in the
sequel. Let $G$ be a connected reductive algebraic group over the algebraic
closure of a finite field $\FF_q$, where $q$ is a prime power.
The variety of unipotent
elements of $G$ is denoted by $G_\uni$.
For a closed subgroup $H$ of $G$ we write
$H^\circ$ for the identity component of $H$ and
$R_\uu(H)$ for the unipotent radical of $H$.
Given $x \in G$  we write
$A(x) = C_G(x)/C_G(x)^\circ$ for the component group of the
centralizer of $x$.

Let $T$ be a maximal torus of $G$ and let $B$ be a Borel subgroup of
$G$ containing $T$. We write $\Psi$ for the root datum of $G$ with
respect to $T$; so $\Psi$ is a quadruple $(X,\Phi,\check X,\check
\Phi)$, where $X$ is the character group of $T$, $\Phi$ is the root
system of $G$ with respect to $T$, $\check X$ is the cocharacter
group of $T$ and $\check \Phi$ is the set of coroots of $G$ with
respect to $T$.  Let $\Pi$ (resp.\ $\Phi^+$) be the base of $\Phi$
(resp.\ the subset of positive roots) determined by $B$. For a
subset $J$ of $\Pi$, the standard parabolic subgroup $P = P_J$ is
defined to be generated by $B$ and the root subgroups $U_{-\alpha}$
for $\alpha \in J$. The unipotent radical of $P$ is denoted by $U =
R_\uu(P)$. We write $L = L_J$ for the unique Levi subgroup of $P$
that contains $T$, i.e.\ $L$ is the standard Levi subgroup of $P$,
see \cite[Prop.\ 1.17]{dignemichel}. The Weyl group of $G$ is
denoted by $W$ and we write $W_L$ for the Weyl group of $L$. As
usual, $l(w)$ denotes the length of $w \in W$ with respect to the
set of Coxeter generators determined by $\Pi$.

For $\beta \in \Phi^+$ write $\beta = \sum_{\alpha \in \Pi}
c_{\alpha,\beta} \alpha$ with $c_{\alpha,\beta} \in \ZZ_{\ge 0}$. We
recall that a prime $p$ is said to be \emph{bad for $\Phi$} if it
divides $c_{\alpha,\beta}$ for some $\alpha$ and $\beta$, else it is
called \emph{good for $\Phi$}; we say $p$ is \emph{good for $G$} if
it is good for $\Phi$.

Now suppose $G$ is defined over $\FF_q$, where $q$ is a power of
$p$, and let $F : G \to G$ be the corresponding Frobenius morphism of $G$.
We write
$G^F$ for the finite group of $F$-fixed points of $G$; likewise for
any $F$-stable subvariety of $G$. Assume that $B$ is $F$-stable
(such a $B$ exists, see \cite[3.15]{dignemichel}) so that $T$ is a
\emph{maximally split} maximal torus of $G$.  For each $w \in W$ we let
$T_w$ be an $F$-stable maximal torus of $G$ obtained by twisting $T$
by $w$, see \cite[3.24]{dignemichel}.

We now recall some notation regarding the representation theory of
the finite group of Lie type $G^F$.  To simplify notation we often
omit the superscript $F$ in places where it should strictly be
written, so for example we write $1_G$ for the trivial character of
$G^F$ and we denote the regular character of $G^F$ by $\chi_G$.
Given an $F$-stable subgroup $H$ of $G$, we write $\Ind_H^G$ for the
induction functor from characters of $H^F$ to characters of $G^F$.
For an $F$-stable Levi subgroup of $G$, we write $R_L^G$ for the
Deligne--Lusztig induction functor from characters of $L^F$ to
characters of $G^F$ (see \cite[\S 11]{dignemichel}).  If $L$ is a
Levi subgroup of an $F$-stable parabolic subgroup $P$ of $G$, then
we recall that $R_L^G$ is the Harish-Chandra induction functor (see
\cite[\S 6]{dignemichel}). For $w \in W$ we write $\Irr(T_w^F)$ for
the character group of $T_w^F$.

In some parts of this paper we consider $l$-adic cohomology groups.
Given a variety $V$ defined over $\FF_q$ we write $H^i(V,\barQl)$
for the $i$th $l$-adic cohomology group of $V$, where $l$ is a prime
not dividing $q$.  We write $H_c^i(V,\barQl)$ for the $l$-adic
cohomology with compact support groups.  Recall that in case $V$ is
projective, the cohomology groups $H^i(V,\barQl)$ and
$H_c^i(V,\barQl)$ coincide.

\subsection{Axiomatic setup for connected reductive algebraic groups}
\label{s:setup} In order to be able to state the results indicated
in the introduction precisely, we use an axiomatic formalism for
connected reductive algebraic groups, following \cite[\S 4]{geck}.  The idea
is that a tuple of combinatorial objects is used to define a family
of connected reductive groups indexed by prime powers. We refer the reader to
\cite[\S 0, \S 3]{dignemichel} for some of the results used below.

Let $\Psi = (X, \Phi, \check X, \check\Phi)$ be a root datum.  Then
given a finite field $\FF_q$, the root datum $\Psi$ determines a
connected reductive algebraic group $G$ over $\bar \FF_q$ and a maximal torus
$T$ of $G$ such that $\Psi$ is the root datum of $G$ with respect to
$T$; moreover, $G$ is unique up to isomorphism and $T$ is unique up
to conjugation in $G$.

Now let $\Pi$ be a base for $\Phi$; this determines a Borel subgroup
$B$ of $G$ containing $T$.  Further, $B$ is determined up to
conjugacy in $G$.

Let $F_0 : X \to X$ be an automorphism of finite order such that
$F_0(\Phi) = \Phi$, $F_0(\Pi) = \Pi$ and $F_0^*(\check\Phi) =
\check\Phi$.  Then for any prime power $q$, the automorphism $F_0$
defines a Frobenius morphism $F: G \to G$ such that the induced
action of $F$ on $X$ is given by $q \cdot F_0$.  Further, $B$ and
$T$ are $F$-stable, so that $T$ is a maximally split maximal torus
of $G$.

Now let $J$ be a subset of $\Pi$ such that $F_0(J) = J$.  Then $J$
determines the standard parabolic subgroup $P = P_J$ of $G$.  If $q$
is a prime power and $F$ is the corresponding Frobenius morphism
$q \cdot F_0$ of $G$, then $P$ is $F$-stable.

Summing up, the discussion above implies that the quadruple $\Delta =
(\Psi,\Pi,F_0,J)$, along with a prime power $q$ determines
\begin{itemize}
\item a connected reductive algebraic group $G$ defined over $\FF_q$ with
corresponding Frobenius morphism $F$ and maximally split $F$-stable
maximal torus $T$;
\item an $F$-stable Borel subgroup $B \supseteq T$ of $G$; and
\item an $F$-stable parabolic subgroup $P \supseteq B$.
\end{itemize}
The connected reductive
algebraic group $G$ is unique up to isomorphism defined over
$\FF_q$; the Borel subgroup $B$ is determined up to conjugacy in
$G^F$; the maximal torus $T$ is determined up to conjugacy in $B^F$;
and (given $B$) $P$ is uniquely determined.

At some points in the sequel we do not require the parabolic
subgroup $P$.  In such situations we only consider the data $\Gamma
= (\Psi,\Pi,F_0)$, which along with a prime power $q$ determines
everything above except for $P$.

We note that if we do not assume that $F_0$ stabilizes a base $\Pi$
of $\Phi$, then the data $(\Psi,F_0)$ and $q$ still determine a
connected reductive algebraic group $G$ that is defined over $\FF_q$
and an $F$-stable maximal torus $T$ of $G$, which is not, in
general, maximally split.  However, there always exists an
$F$-stable Borel subgroup of $G$ (see \cite[3.15]{dignemichel}).
Therefore, it does no harm to make the assumption that $F_0$
stabilizes $\Pi$.

The notation we use for $G$, $B$, $T$ and $P$ does not reflect the
fact that their $\FF_q$-structure depends on the choice of a prime
power $q$. Let $q$ be a prime power and $m$ a positive integer,
write $F$ for the Frobenius morphism corresponding to $q$.  Then it
is not necessarily the case that the Frobenius morphism
corresponding to the prime power $q^m$ is $F^m$, i.e.\ the
definition of $G$ over $\FF_{q^m}$ is not necessarily obtained from
the $\FF_q$-structure by extending scalars. For example, suppose $G
= \SL_n$ with definition over $\FF_q$ so that $G^F$ is the special
unitary group $\SU_n(q^2)$, i.e.\ $F_0$ is given by the graph
automorphism of the Dynkin diagram of type $A_{n-1}$. Then the
automorphism $F_0^2$ is the identity, so that the group of
$\FF_{q^2}$-rational points of $G$ is $G^{F^2} = \SL_n(q^2)$ and not
$\SU_n(q^4)$. However, in order to keep the notation short, we
choose not to show this dependence on $q$.

\begin{rem} \label{R:varyq}
Despite the lack of outward exhibition of dependence on $q$ discussed above,
in the
sequel we wish to allow $q$ to vary and for $G$, $B$, $T$ and $P$ to
define groups and $F$ a Frobenius morphism of $G$ for each $q$. For
example, this convention is required in order to make sense of statements
such as the following: ``the number of $F$-stable unipotent
conjugacy classes of $G$ is independent of $q$'', or ``the number of conjugates
of $P^F$ having non-trivial intersection with a given $F$-stable
conjugacy class of unipotent elements is given by a polynomial in
$q$''.
\end{rem}

Now we discuss some further information that is determined by the
data $\Delta = (\Psi,\Pi,F_0,J)$ and prime power $q$, where $\Psi =
(X, \Phi, \check X, \check\Phi)$ is a root datum. First we note that
the unipotent radical $U = R_\uu(P)$ of the $F$-stable parabolic
subgroup $P = P_J$ of $G$ and the unique Levi subgroup $L = L_J$ of
$P$ containing $T$ are determined. Since both $P$ and $T$ are
$F$-stable, so is $L$.
The Weyl group $W$ of $G$ and the Weyl group $W_L$ of $L$ are
determined by $\Delta$; further, we note that $W$ and $W_L$ are
independent of $q$.  The $F$-stable twisted tori $T_w$ ($w \in W$)
are determined up to conjugacy in $G^F$.  The root lattice of
$\Phi$, which is denoted by $\ZZ\Phi$, is determined by $\Delta$.

Lastly, in this section we discuss two further properties of the
group $G$ that are determined by $\Delta$.   At some points in the
sequel we wish to assume that the centre of $G$ is connected.  We
note that this being the case for all $q$ corresponds to $X/\ZZ\Phi$
being torsion free; this is because the character group of $Z(G)$ is
$X/\ZZ\Phi$ modulo its $p$-torsion subgroup, where $q$ is a power of
the prime $p$. Also at some points we assume that $G$ is simple
modulo its centre, we note that this corresponds to $\Phi$ being an
irreducible root system.

\subsection{On the unipotent conjugacy classes of $G^F$} \label{s:classes}

In this subsection we use the axiomatic setup from \S \ref{s:setup}.
In particular, we fix a triple $\Gamma = (\Psi,\Pi,F_0)$, where
$\Psi = (X,\Phi,\check X,\check \Phi)$, so that for
each choice of a prime power $q$ we obtain a connected
reductive group $G$ defined over $\FF_q$; we write $F : G \to G$ for
the corresponding Frobenius morphism of $G$. We use the convention
discussed in Remark \ref{R:varyq} to vary $q$.
Further, any prime power that we
consider is assumed to be a power of a good prime for $\Phi$.

It is known that the parameterization of the unipotent classes of $G$
and the structure of the associated component groups $A(u) =
C_G(u)/C_G(u)^\circ$ for $u \in G_\uni$ is independent of $q$ (cf.\
\cite{pomerening}, \cite{premet}, \cite{mcnsom}).  We would like a
similar statement for the $F$-stable unipotent classes of $G$, and
for the unipotent $G^F$-classes.

There is a parameterization of the $F$-stable unipotent classes as
follows (see \cite{shoji}, \cite{kawanaka}, and also \cite[\S
4]{geck}). There is a finite set $E$ depending only on $\Gamma$ (and
not on $q$) so that there is a map $E \to G_\uni^F$ written $e
\mapsto u_e$ so that $\{u_e \mid e \in E\}$ is a complete set of
representatives for the $F$-stable unipotent classes in $G$.
Therefore, the $F$-stable conjugacy classes are parameterized by $E$
and this parameterization is independent of $q$.

In order to get a parameterization of the unipotent $G^F$-conjugacy classes
we require some additional assumptions.  We
assume that $X/\ZZ\Phi$ is torsion free, so that the
centre of $G$ is connected, and $\Phi$ is irreducible, so that $G$
is simple modulo its centre.

With these additional assumptions, there is a collection of groups
$\{A_e \mid e \in E\}$ (not depending on $q$) such that the
component group $A(u_e)$ is isomorphic to $A_e$ for every $e \in E$
and the induced action of $F$ on $A(u_e)$ is trivial. Consequently,
for $e \in E$, the $G^F$-conjugacy classes in $(G \cdot  u_e)^F$ are
in correspondence with the conjugacy classes of $A_e$, cf.\ \cite[I,
2.7]{springersteinberg}. So the unipotent $G^F$-classes are
parameterized by pairs $(e,c)$, where $e \in E$ and $c$ is a
conjugacy class of $A_e$.  Moreover, this parameterization is
independent of $q$.

These results largely follow from work of T.~Shoji on split
elements. Following \cite[\S5]{shoji}, an element $u \in G_\uni^F$
is called \emph{split} provided each irreducible component of
$\CB_u$ (see \S \ref{s:purity} for a definition of $\CB_u$) is
stable under the action of $F$. If $q$ is a power of a good prime
for $\Phi$, then split elements exist in each $F$-stable unipotent
$G$-conjugacy class and are unique up to $G^F$-conjugacy with one
exception, see \cite[Rem.\ 5.1]{shoji}. If $G$ is of type $E_8$, $q
\equiv -1 \mod 3$, and $C$ is the class $E_8(b_6)$ in the
Bala--Carter labeling (this class is $D_8(a_3)$ in Mizuno's labeling
\cite{mizuno}), then there are no split elements in $C^F$. If $u$ is
split, then $F$ acts trivially on $A(u)$, see \cite{shoji}.  The one
case for $G$ of type $E_8$ where there is no split element can be
dealt with explicitly, see \cite{kawanaka}.

Thus we have the following result, where we are using the convention
discussed in Remark \ref{R:varyq} to allow $q$ to vary. This
proposition is also stated in \cite[\S 4.1]{geck}.

\begin{prop}
\label{prop:classes} Fix the data $\Gamma = (\Psi,\Pi,F_0)$, where
$\Psi = (X, \Phi, \check X, \check\Phi)$. Assume that $\Phi$ is
irreducible and $X/\ZZ\Phi$ is torsion free.  For a prime power $q$,
let $G$ be the connected reductive algebraic group and $F$ the Frobenius
morphism determined by $\Gamma$. Then the parameterization of the
unipotent conjugacy classes in $G^F$ is independent of $q$, for $q$
a power of a good prime for $\Phi$.
\end{prop}

We now turn our attention to showing that Proposition
\ref{prop:classes} remains true when we remove the assumption
that $G$ is simple modulo its centre, i.e.\
that the root system $\Phi$ is possibly reducible.  At the end of this
subsection we briefly discuss what happens when $G$
has disconnected centre.

We begin by assuming that $G$ is semisimple and of adjoint type,
equivalently that $X = \ZZ\Phi$. In this case the root system $\Phi$
has irreducible components on which the automorphism $F_0$ acts. We
may write $\Phi$ as a disjoint union $\Phi = \Phi_1 \cup \dots \cup
\Phi_s$ of $F_0$-stable root systems, where the irreducible
components of each $\Phi_i$ form a single orbit under the action of
$F_0$. This gives rise to a direct product decomposition $G = G_1
\times \dots \times G_s$, where each $G_i$ is a semisimple algebraic
group of adjoint type with root system $\Phi_i$; moreover each $G_i$
is $F$-stable. Now for fixed $i$, the root system $\Phi_i$ is the
disjoint union of irreducible root systems of the same type; let
$\Upsilon_i$ be a root system of this type.  Then $G_i$ is a direct
product $G_i = H_i \times \dots \times H_i$ where $H_i$ is a simple
algebraic group of adjoint type with root system $\Upsilon_i$.  Let
$l_i$ be the number of factors $H_i$ in the direct product
decomposition of $G_i$.  Then we see that $G_i^F \cong
H_i^{F^{l_i}}$, so that the $\FF_q$-rational points of $G_i$ are in
correspondence with the $\FF_{q^{l_i}}$-rational points of $H_i$. It
follows that the unipotent $G_i^F$-conjugacy classes have the same
parameterization as the unipotent $H_i^{F^{l_i}}$-conjugacy classes,
i.e.\ there is a finite set $E_i$ and for each $e \in E_i$ there is
a group $A_{i,e}$ such that the unipotent $G_i^F$-conjugacy classes
are given by pairs $(e,c)$ where $e \in E_i$ and $c$ is a conjugacy
class of $A_{i,e}$. More precisely, there exists $u_{i,e} \in G_i$
for $e \in E_i$ such that the $F$-stable $G_i$-conjugacy classes are
precisely the classes $G_i \cdot u_{i,e}$ and the component group
$A(u_{i,e})$ of the centralizer of $u_{i,e}$ is $A_{i,e}$.  Let $E =
E_1 \times \dots \times E_s$ and for $e = (e_1,\dots,e_s) \in E$ let
$A_e = A_{1,e_1} \times \dots \times A_{s,e_s}$.  Then we see that
the unipotent $G^F$-conjugacy classes are parameterized by pairs
$(e,c)$, where $e \in E$ and $c$ is a conjugacy class of $A_e$.  For
$e \in E$, we let $u_e$ be a representative of the $F$-stable
unipotent $G$-class, such that $A(u_e) \cong A_e$ and $F$ acts trivially
on $A(u_e)$.

We now consider the case where the centre of $G$ is not necessarily
trivial, but is connected.  In this case we let $\hat G$ be the
semisimple adjoint group of the same type as $G$ and let $\phi : G
\to \hat G$ be the natural map. The Frobenius morphism $F$ of $G$
induces a Frobenius morphism of $\hat G$, which we also denote by
$F$. We note that $\phi$ induces a bijective morphism $G_\uni \to
\hat G_\uni$, and this induces a bijection between the $F$-stable
unipotent conjugacy classes of $G$ and those of $\hat G$. Let $u \in
G_\uni^F$ be a representative of an $F$-stable unipotent
$G$-conjugacy class so that $\phi(u) = u_e$ for some $e \in E$ as
above. We have the exact sequence of algebraic groups
\begin{equation} \label{eqses}
\{1\} \to Z(G) \to C_G(u) \to C_{\hat G}(\phi(u)) \to \{1\}.
\end{equation}
Further we note that $F$ commutes with the maps in this sequence.
Since taking component groups defines a right exact functor on the
category of algebraic groups over $\bar \FF_q $ and $Z(G)$ is
connected, we get an exact sequence of finite groups
\[
\{1\} \to A(u) \to \hat A(\phi(u)) \to \{1\}.
\]
Moreover, $F$ commutes with this isomorphism.  Therefore, it follows
that $\phi$ induces a bijection between the unipotent
$G^F$-conjugacy classes and the unipotent $\hat G^F$-conjugacy
classes.

The discussion above gives the following generalization of
Proposition \ref{prop:classes}.

\begin{prop}
\label{prop:classes2} Fix the data $\Gamma = (\Psi,\Pi,F_0)$, where
$\Psi = (X, \Phi, \check X, \check\Phi)$. Assume that $X/\ZZ\Phi$ is
torsion free. For a prime power $q$, let $G$ be the reductive
algebraic group and $F$ the Frobenius morphism determined by
$\Gamma$. Then the parameterization of the unipotent conjugacy
classes in $G^F$ is independent of $q$, for $q$ a power of a good
prime for $\Phi$.
\end{prop}

\begin{rem} \label{R:varyq2}
The discussion proving Proposition \ref{prop:classes2} gives a
parameterization of the unipotent $G^F$-conjugacy classes given by
pairs $(e,c)$, where $e \in E$ and $c$ is a conjugacy class of
$A_e$.  We may let $u = u_{e,c} \in G_\uni^F$ be a representative of
this conjugacy class.  We now use the convention in Remark
\ref{R:varyq} to vary $q$.  In general, as we vary $q$, the
representative $u$ changes.  For the statement of results in
Sections \ref{S:ratpts} and \ref{s:alperin}, we use the convention
that when we fix a representative of a unipotent class and let $q$
vary, then it is understood that $u$ also changes.
\end{rem}

To finish this subsection we briefly consider the situation when $G$ has
disconnected centre. As we do not require this case in the sequel,
the discussion below is brief and omits details.

Suppose that $G$ has disconnected centre.
Let $\hat G$ be the semisimple adjoint group
of the same type as $G$ and let $\phi : G \to \hat G$ be the natural
map.  We still have a bijective morphism $G_\uni \to \hat G_\uni$
induced by $\phi$ and this affords a bijection between the $F$-stable
unipotent classes in $G$ and $\hat G$.  The
arguments above using the exact sequence \eqref{eqses} imply that
taking the quotient of $G$ by $Z(G)^\circ$ does not affect the
unipotent $G^F$-conjugacy classes.  Therefore, we may assume that
$Z(G)$ is finite.

Let $u \in G_\uni^F$ be such that $\phi(u) = u_e$ for some $e \in E$
as above.  We still have the
short exact sequence \eqref{eqses}.
When we take component groups, we get the exact sequence
\begin{equation*}
Z(G) \to A(u) \to \hat A(\phi(u)) \to \{1\},
\end{equation*}
and the action of $F$ commutes with the maps in this exact sequence.
Let $Z$ be the image of $Z(G)$ in $A(u)$.  Then we have $A(u)/Z \cong
\hat A(\phi(u))$ and $Z$ is stable under $F$-conjugation by $A(u)$.
It follows
that the $F$-conjugacy classes of $A(u)$ are given by pairs $(c,d)$, where
$c$ is a conjugacy class of $\hat A(\phi(u))$ with representative
$xZ \in A(u)/Z \cong \hat A(\phi(u))$, and
$d$ is an $F$-conjugacy class of $C_{A(u)}(x)$ in $Z$.
This in turn says that the $G^F$-conjugacy classes in $(G \cdot u)^F$
are parameterized by pairs $(c,d)$ as above.
In general, this parameterization depends on $q$ but only up to
congruences of $q$ modulo some positive integer.

We illustrate this splitting of conjugacy classes
in the disconnected centre case with a simple example.

\begin{exmp}
\label{ex:SL3}
Suppose $G = \SL_3$ with the standard definition over $\FF_q$.
Then $G$ has three unipotent conjugacy classes labeled by the
partitions $(3)$, $(2,1)$ and $(1,1,1)$ of $3$.  The regular unipotent class
corresponds to the partition $(3)$ and a split representative is
\[
u = \left(\begin{array}{ccc}
1 & 1 & 0 \\ 0 & 1 & 1 \\ 0 & 0 & 1
\end{array} \right).
\]
One can easily calculate the centralizer of $u$, and deduce that
$A(u)$ is the group of cube roots of unity in $\bar \FF_q$ and the
action of $F$ on $A(u)$ is the $q$th power map.  If $q$ is congruent
to $1$ modulo $3$, then the action of $F$ on $A(u)$ is trivial and
so there are three $F$-conjugacy classes in $A(u)$.  If $q$ is
congruent to $-1$ modulo $3$, then there is a single $F$-conjugacy
class in $A(u)$ and if  $q$ is a power of $3$, then $A(u)$ is the
trivial group and thus is a single $F$-conjugacy class. It follows
that the $G^F$-conjugacy classes in $(G \cdot u)^F$ depend on the
congruence class of $q$ modulo $3$.
\end{exmp}

\subsection{On purity of generalized flag varieties}
\label{s:purity} In this subsection we use the notation from \S
\ref{s:notn}; we do not require the axiomatic formalism from \S
\ref{s:setup}.

Let $\CB$ be the flag variety of $G$ consisting of all Borel
subgroups of $G$.  Fix a parabolic subgroup $P$ of $G$.
Let $\CP = \{{}^gP \mid g \in G\}$
be the corresponding \emph{generalized flag variety} consisting of
all the $G$-conjugates of $P$.

For a unipotent element $u \in G$, we consider the fixed point
subvarieties $\CB_u = \{B' \in \CB \mid u \in B'\}$ of $\CB$ and
$\CP_u = \{P' \in \CP \mid u \in P'\}$ of $\CP$. We also consider
the related subvariety $\CP_u^0 = \{P' \in \CP \mid u \in
R_\uu(P')\}$ of $\CP_u$. Note that for $P = B$ we have $\CP_u =
\CP_u^0 = \CB_u$. We can identify the variety $\CB_u$ with a fibre
of Springer's resolution of the singularities of the unipotent
variety $G_\uni$ of $G$. The fixed point varieties $\CP_u$ can be
viewed as fibres of \emph{partial resolutions} of $G_\uni$, and have
been studied by Borho--MacPherson \cite{borhomac} and Spaltenstein
\cite{spaltenstein}. We note that while the results in
\cite{borhomac} only apply in characteristic $0$, Spaltenstein's
treatment \cite{spaltenstein} is the positive characteristic setting
we require.

Let $V$ be a variety defined over $\FF_q$ with Frobenius morphism
$F$ corresponding to the $\FF_q$-structure. We recall that $V$ is
called  \emph{pure} (or is said to satisfy the \emph{purity
condition}), if the eigenvalues of $F$ on the $l$-adic cohomology
groups $H_c^i(V,\barQl)$ are algebraic integers of absolute value
$q^{i/2}$.

Let $u \in G_\uni^F$. In \cite{springer1}, Springer showed that
despite the fact that the fixed point varieties $\CB_u$ are not
smooth, they are pure. We require this purity condition also for the
varieties $\CP_u$ and $\CP_u^0$.

Let $L$ be the unique Levi subgroup of $P$ containing $T$.
It follows from a result of Borho--MacPherson \cite[Cor.\
2.8]{borhomac} (see also \cite[Thm.\ 1.2]{spaltenstein}) that
for each $i$ there is an isomorphism of
$\barQl$-spaces
\begin{equation}
\label{e:pu}
H^i(\CP_u, \barQl) \cong H^i(\CB_u, \barQl)^{W_L},
\end{equation}
where $H^i(\CB_u, \barQl)^{W_L}$ is the subspace of $H^i(\CB_u,
\barQl)$ of $W_L$-invariants (where the action of $W$ on $H^i(\CB_u,
\barQl)$ is given by the action of $W$ defined by Springer
\cite{springer1} multiplied by the sign character of $W$). The
isomorphism in \eqref{e:pu} is induced by the natural map $\CB_u \to
\CP_u$ (cf.\ \cite[Thm.\ 1.2]{spaltenstein}), which is clearly
$F$-equivariant. Therefore, the isomorphism in \eqref{e:pu} is
compatible with the induced $F$-actions on both sides. As a
consequence, the eigenvalues of $F$ on $H^i(\CP_u, \barQl)$ are
algebraic integers of modulus $q^{i/2}$ so that $\CP_u$ is pure.

Since $P$ and $T$ are $F$-stable, so is $L$ and the
induced action of $F$ on $L$
is a Frobenius morphism. Whence $F$ acts on $\CB(L)$, the
variety of Borel subgroups of $L$. Let $d_0 = \dim \CB(L)$.
Since $\CB(L)$ is smooth, and
$H^{2d_0}(\CB(L),\barQl) \cong \barQl$,
the eigenvalue of $F$ on this
one-dimensional space is $q^{d_0}$.
Thanks to \cite[Thm.\  3.3, Cor.\ 3.4(b)]{borhomac},
we have an isomorphism of
$\barQl$-spaces
\begin{equation}
\label{e:pu0}
H^{i-2d_0}(\CP_u^0, \barQl) \otimes H^{2d_0}(\CB(L), \barQl)
\cong H^i(\CB_u, \barQl)^{\epsilon W_L},
\end{equation}
for each $i$, where $H^i(\CB_u, \barQl)^{\epsilon W_L}$ is the
subspace of $H^i(\CB_u, \barQl)$ on which $W_L$ acts via the sign
representation $\epsilon$ (i.e. the subspace of
$W_L$-anti-invariants). Further, the isomorphism in \eqref{e:pu0} is
$F$-equivariant. By Springer's purity result, the eigenvalues of $F$
on $H^i(\CB_u, \barQl)^{\epsilon W_L}$ are algebraic integers of
modulus $q^{i/2}$. Consequently, the eigenvalues of $F$ on
$H^{i-2d_0}(\CP_u^0, \barQl)$ are algebraic integers of modulus
$q^{(i-2d_0)/2}$ for all $i$. It follows that $\CP_u^0$ is pure.

\subsection{Deligne--Lusztig generalized characters and Green functions}
\label{s:dlgenchars}

For this subsection we use the axiomatic formalism given in \S
\ref{s:setup}. We fix the triple $\Gamma = (\Psi,\Pi,F_0)$, where
$\Psi = (X, \Phi, \check X, \check\Phi)$, so for each prime power
$q$ we have a connected reductive algebraic group $G$ defined over
$\FF_q$ with Frobenius morphism $F$.  Further, we have a maximal
torus $T$ of $G$ and the Weyl group $W$ of $G$ with respect to $T$.
We recall from \S \ref{s:notn} that for a twisted torus $T_w$ ($w
\in W$), we write $R_{T_w}^G$ for the Deligne--Lusztig induction
functor from characters of $T_w^F$ to characters of $G^F$.  In this
subsection we recall some known results about the values of
$R_{T_w}^G(\theta)$ on unipotent elements of $G^F$, where $\theta
\in \Irr(T_w^F)$.

Let $w \in W$ and $\theta \in \Irr(T_w^F)$.  We begin by recalling
that as a function $G_\uni^F \to \bar\QQ$ the restriction
$R_{T_w}^G(\theta)$ to $G_\uni^F$ does not depend on $\theta$, see
\cite[Cor.\ 7.2.9]{carter}. The restriction of $R_{T_w}^G(\theta)$
to $G_\uni^F$ is called a \emph{Green function} and is denoted by
$Q_{T_w}^G$, see \cite[\S 7.6]{carter}.

Suppose $q$ is a power of a good prime for $\Phi$ and assume that
$X/\ZZ\Phi$ is torsion free (i.e.\  the centre of $G$ is connected),
so that Proposition \ref{prop:classes2} applies.  With these
assumptions, work of Shoji \cite[Prop.\ 6.1]{shoji} tells us that
the Green functions $Q_{T_w}^G$ evaluated on a unipotent class of
$G^F$ are given by certain polynomials in $q$ whose coefficients are
independent of $q$, if $G$ does not have a simple component of type
$E_8$. If $G$ has a simple component of type $E_8$, then for $u \in
G_\uni^F$, the value of $Q_{T_w}^G(u)$ is given by one of two
polynomials in $q$ depending on whether $q$ is congruent $1$ or $-1$
modulo $3$, see \cite[Rem.\ 6.2]{shoji}.

\begin{rem}
\label{r:funnyE8} For $G$ simple of type $E_8$, we only require two
polynomials in case $u$ lies in the conjugacy class with
Bala--Carter label $E_8(b_6)$, see \cite[Rem.\ 5.1]{shoji}.  This is
precisely the conjugacy class for which there is no split element if
$q$ is congruent to $-1$ mod $q$, as mentioned in \S
\ref{s:classes}.  For reductive $G$ the two polynomials are only
required in case there is a simple component $H$ of $G$ of type
$E_8$, such that the size of the $F$-orbit of $H$ is odd and the
projection of $u$ into $H$ belongs to the class with Bala--Carter
label $E_8(b_6)$. For, if the $F$-orbit of $H$ is of even size $2m$,
then the group of $F$-stable points of the direct product of the
simple groups in this orbit is a finite simple group of type $E_8$
over the field $\FF_{q^{2m}}$, so that there does exist a split
element (since $q^{2m} \equiv 1 \mod 3$).
\end{rem}

\begin{rem} \label{rem:badchargreen}
In case $q$ is a power of a bad prime for $G$, it is conjectured
that the values of the Green functions $Q_{T_w}^G(u)$ are also given
by polynomials in $q$; this has been verified in many cases.  In
general, these polynomials differ from those for good
characteristic. Based on results of Shoji (\cite{shoji}), it is
shown by M.~Geck in \cite[Prop.\ 3.5]{geck} that the algorithm for
the computation of Green functions presented by G.~Lusztig in
\cite[Chap.\ 24]{LusztigV}, works in any characteristic. In order to
show that the Green functions are polynomials in bad characteristic,
one has to find suitable analogues of the split representatives
which exist in good  characteristic. The latter problem is open in
general.
\end{rem}

\subsection{Two standard lemmas} \label{s:twolemmas}

In this subsection we recall two standard lemmas that we require in the sequel.
The first of these concerns a condition on a rational function
that forces it to be a
polynomial, we include a proof of this well-known result for the reader's
convenience.

\begin{lem} \label{L:rattopoly}
Let $f(z) \in \QQ(z)$ be a rational function.  Suppose that there are
infinitely many $n \in \ZZ$ such that $f(n) \in \ZZ$.  Then $f(z) \in \QQ[z]$
is in fact a polynomial.
\end{lem}

\begin{proof}
We may write $f(z) = g(z) + \frac{h(z)}{k(z)}$, where $g(z),h(z),k(z)
\in \QQ[z]$ and the degree of $h(z)$ is strictly less than that of $k(z)$.
Suppose for a contradiction that $h(z) \ne 0$.
We can find a positive integer $N$ such that $Ng(z) \in \ZZ[z]$.  Since
$f(n) \in \ZZ$ for infinitely many $n \in \ZZ$, there exists such $n$ for which
$0 < \frac{h(n)}{k(n)} < 1/N$.  But then we have $Nf(n) =
Ng(n) + N\frac{h(n)}{k(n)}$, and the left hand side is an integer, whereas
the right hand side cannot be an integer.  This contradiction shows that
$h(z) = 0$, so that $f(z) = g(z) \in \QQ[z]$, as required.
\end{proof}

For the next lemma we refer the reader to \S \ref{s:purity} for the
definition of a pure variety. Using the Grothendieck trace formula
(e.g.\ see \cite[Thm.\ 10.4]{dignemichel}), Lemma \ref{L:intcoeffs}
can be proved with standard arguments, see for example \cite[\S
5]{beyspalt}.

\begin{lem}
\label{L:intcoeffs}
Let $V$ be a pure variety defined over
a finite field $\FF_q$. Let $F$ be the Frobenius morphism
corresponding to the $\FF_q$-structure on $V$.
Suppose there exists $g(z) \in \QQ[z]$ such that $|V^{F^s}| = g(q^s)$ for
any $s \in \ZZ_{\ge 1}$.  Then
\[
g(z) = \sum_{i=0}^{\dim V} \dim H_c^{2i}(V,\barQl) z^i.
\]
In particular, $g(z) \in \ZZ[z]$ and the coefficients in $g(z)$ are
positive integers. Further, we have $H_c^i(V,\barQl) = \{0\}$ for
$i$ odd, and all the eigenvalues of $F^s$ on $H_c^{2i}(V,\barQl)$
are $q^{is}$ for all $i = 0, \ldots, \dim V$ and all $s \in \ZZ_{\ge
1}$.
\end{lem}

\section{Rational points of fixed point subvarieties of
generalized flag varieties}
\label{S:ratpts}

In this section we prove Theorem \ref{thmmain} and Proposition
\ref{p:betti}.  The proofs of these results are contained in \S
\ref{s:poly} and \S \ref{s:betti}.  Theorem \ref{thmmain} is a
consequence of Proposition \ref{prop1}, which is proved by combining
well-known results from the representation theory of finite groups
of Lie type.

\subsection{Fixed point subvarieties of flag varieties}
\label{Sfixedpoint}  In this subsection we use the notation given in
\S \ref{s:notn} and also the notation from \S \ref{s:purity}; we do
not require the axiomatic setup from \S \ref{s:setup}. In
particular, $B$ is an $F$-stable Borel subgroup of $G$ containing an
$F$-stable maximal torus $T$, and $P$ is an $F$-stable parabolic
subgroup of $G$ containing $B$.  The unipotent radical of $P$ is
denoted by $U$. Note that $U$ is also $F$-stable.

The Frobenius morphism $F$ acts on the variety $\CP$ of conjugates
of $P$. For $u \in G_\uni^F$ this induces an action of $F$ on the
fixed point varieties $\CP_u$ and also on $\CP_u^0$.  Therefore, we
may consider the $F$-fixed points of each of these varieties:
$\CP^F$, $\CP_u^F$ and $(\CP_u^0)^F$.  We may identify $\CP$ with
$G/P$, because $P = N_G(P)$.  Then by \cite[Cor.\ 3.13]{dignemichel}
we identify $\CP^F$ with $G^F/P^F$, which in turn we can identify
with $\{{}^gP^F \mid g \in G^F\}$, because $N_{G^F}(P^F) = P^F$.
This allows us to identify $\CP_u^F$ with $\{{}^gP^F \mid g \in G^F,
u \in {}^gP^F\}$. In a similar way, we identify $(\CP_u^0)^F$ with
$\{{}^gP^F \mid g \in G^F, u \in {}^gU^F\}$ which in turn can be
identified with $\{{}^gU^F \mid g \in G^F, u \in {}^gU^F\}$.

For an $F$-stable subgroup $H$ of $G$ and for $x \in G^F$ we define
\[
f_H^G(x) = |\{{}^gH^F \mid g \in G^F, x \in {}^gH^F\}|.
\]
To keep the notation short, we do not show the dependence on $F$
in the notation $f_H^G(x)$.
Note that if $x \in G^F$ is not conjugate to an element in $H^F$,
then $f_H^G(x) = 0$.

The following lemma is a consequence of the discussion above.

\begin{lem}
\label{lem3}
Let $u \in G^F_\uni$. Then we have
\begin{itemize}
\item[(i)] $|\CP_u^F| = f_P^G(u)$;
\item[(ii)] $|(\CP_u^0)^F| = f_U^G(u)$.
\end{itemize}
\end{lem}

Let $L$ be the unique Levi subgroup of $P$ containing $T$.  Since
$P$ and $T$ are $F$-stable, so is $L$.  For the statements of the
next lemmas we recall that $W_L$ denotes the Weyl group of $L$ and
$T_w$ is a maximal torus obtained from $T$ by twisting by $w \in W$.

\begin{lem}
\label{lem4} Let $u \in G^F_\uni$. Then we have
\[
f_P^G(u) = \frac{1}{|W_L|} \sum_{w \in W_L}Q_{T_w}^G(u).
\]
\end{lem}

\begin{proof}
First note that  $f_P^G(u)$ is the value at $u$ of the permutation
character $\Ind_P^G(1_P)$ of $G^F$ on the cosets of $P^F$. Now
$\Ind_P^G(1_P)$ is the trivial character $1_L$ of $L^F$
Harish-Chandra induced to $G^F$, i.e.\ as a function $G_\uni^F \to
\bar\QQ$ we have
\[
f_P^G  = R_L^G(1_L),
\]
where $R_L^G$ is the Harish-Chandra induction functor. By
\cite[Prop.\ 12.13]{dignemichel} the trivial character $1_L$ of
$L^F$ is given in terms of Deligne--Lusztig generalized characters
as follows
\[
1_L = \frac{1}{|W_L|} \sum_{w \in W_L} R_{T_w}^L(1_{T_w}),
\]
where $R_{T_w}^L$ is the Deligne--Lusztig induction functor.
Applying the Harish-Chandra induction functor $R_L^G$ to
$R_{T_w}^L(1_{T_w})$ gives the corresponding generalized character
$R_{T_w}^G(1_{T_w})$ for $G^F$, cf.\ \cite[Prop.\ 4.7]{dignemichel}.
Therefore, we have
\[
R_L^G(1_L) = \frac{1}{|W_L|} \sum_{w \in W_L} R_{T_w}^G(1_{T_w}).
\]
As discussed in \S \ref{s:dlgenchars} the restriction of the virtual
character $R_{T_w}^G(1_{T_w})$ to the unipotent elements of
$G^F$ is just the Green function $Q_{T_w}^G$ of $G^F$.  Thus we obtain
\[
f_P^G(u) = \frac{1}{|W_L|} \sum_{w \in W_L} Q_{T_w}^G(u),
\]
as claimed.
\end{proof}

The formula in Lemma \ref{lem4} can also be found in \cite[Lem.\ 2.3]{LLS}.
Our next result is proved in a similar way to Lemma \ref{lem4}.

\begin{lem}
\label{lem5} Let $u \in G^F_\uni$. Then we have
\[
f_U^G(u) = \frac{1}{|L^F|_{p}|W_L|}\sum_{w \in W_L}
(-1)^{l(w)} Q_{T_w}^G(u).
\]
\end{lem}

\begin{proof}
First note that $|L^F| \cdot f_{U}^G(u)$ is the value at $u$ of the
permutation character $\Ind_U^G(1_U)$ of $G^F$ on the cosets of
$U^F$, because $N_{G^F}(U^F) = P^F = L^FU^F$. Since $U^F$ is normal
in $P^F$, we see that $\Ind_U^G(1_U)$ is the regular character
$\chi_L$ of $L^F$ Harish-Chandra induced to $G^F$, i.e.
\[
|L^F| \cdot f_U^G  = R_L^G(\chi_L).
\]
By \cite[Cor.\ 12.14]{dignemichel}
the regular character $\chi_L$ of $L^F$ is given by
\begin{equation*}
\chi_{L} = \frac{1}{|W_L|} \sum_{w \in W_L} \dim
(R_{T_w}^L(1_{T_w})) R_{T_w}^L(\chi_{T_w}).
\end{equation*}
Using \cite[Cor.\ 12.9, Rem.\ 12.10]{dignemichel}, we have
\begin{equation*}
\chi_{L} = |L^F|_{p'} \cdot \frac{1}{|W_L|} \sum_{w \in W_L}
(-1)^{l(w)} |T_w^F|^{-1} R_{T_w}^L(\chi_{T_w}).
\end{equation*}
Since $\chi_{T_w} = \sum_{\theta \in \Irr(T_w^F)} \theta$, we have
\begin{equation*}
\chi_{L} = |L^F|_{p'} \cdot \frac{1}{|W_L|} \sum_{w \in W_L}
(-1)^{l(w)} |T_w^F|^{-1} \sum_{\theta \in \Irr(T_w^F)}
R_{T_w}^L(\theta).
\end{equation*}
Thus, the regular character of $L^F$ is a $\mathbb Q$-linear
combination of generalized Deligne--Lusztig characters,
$R_{T_w}^L(\theta)$. Harish-Chandra induction applied to
$R_{T_w}^L(\theta)$ gives the corresponding generalized character
$R_{T_w}^G(\theta)$ for $G^F$.  Therefore, we have
\begin{equation*}
R_L^G(\chi_{L}) = |L^F|_{p'} \cdot  \frac{1}{|W_L|} \sum_{w \in W_L}
(-1)^{l(w)} |T_w^F|^{-1} \sum_{\theta \in \Irr(T_w^F)}
R_{T_w}^G(\theta).
\end{equation*}
On unipotent elements of $G^F$ the virtual character
$R_{T_w}^G(\theta)$ is independent of the choice of the character
$\theta$ of $T_w^F$ and is just the Green function $Q_{T_w}^G$, as
discussed in \S \ref{s:dlgenchars}. Therefore, we have that
$R_{T_w}^G(\theta)(u) = R_{T_w}^G(1_{T_w})(u) = Q_{T_w}^G(u)$. Thus we
obtain
\begin{align*}
f_{U}^G(u) & =
|L^F|_{p'} \cdot \frac{1}{|L^F|} \cdot  \frac{1}{|W_L|} \sum_{w \in W_L}
(-1)^{l(w)}  |T_w^F|^{-1}
|\Irr(T_w^F)| Q_{T_w}^G(u)\\
& = \frac{1}{|L^F|_{p}} \cdot  \frac{1}{|W_L|}\sum_{w \in W_L}
(-1)^{l(w)} Q_{T_w}^G(u),
\end{align*}
as claimed.
\end{proof}

Combining Lemmas  \ref{lem3}--\ref{lem5} we get the following
proposition.

\begin{prop}
\label{prop1}
Let $P$ be an $F$-stable parabolic subgroup of $G$
containing an $F$-stable maximal torus $T$ of $G$. Let $L$ be the
Levi subgroup of $P$ containing $T$ and let $u \in G^F_\uni$. Then
\begin{itemize}
\item[(i)]
\[
|\CP_u^F| = \frac{1}{|W_L|} \sum_{w \in W_L}Q_{T_w}^G(u); 
\]
\item[(ii)]
\[
|(\CP_u^0)^F| = \frac{1}{|L^F|_{p}|W_L|} \sum_{w \in W_L}
(-1)^{l(w)} Q_{T_w}^G(u).
\]
\end{itemize}
\end{prop}

We note that the special case $\CP = \CB$ of Proposition \ref{prop1}
is well-known: $|\CB_u^F| = Q_T^G(u) = R_T^G(1_T)(u)$, e.g.\
see \cite[\S 5]{springer1}.

We recall that two parabolic subgroups of $G$ are called
\emph{associated} if they have Levi subgroups that are conjugate in
$G$.  The formulas for $|\CP_u^F|$ and $|(\CP_u^0)^F|$ given in
Proposition \ref{prop1} depend only on the Levi subgroup $L$ and not
on the parabolic subgroup containing $L$; this is due to the fact
that the Harish-Chandra induction functor $R_L^G$ does not depend on
the choice of parabolic subgroup containing $L$ that is used to
define it, \cite[Prop.\ 6.1]{dignemichel}.  Therefore, we deduce the
following corollary; we note that this was observed also
by G.~Lusztig, see \cite[II 4.16]{spaltensteinbook}.

\begin{cor}
\label{cor6}
Let $P$ and $Q$ be associated $F$-stable parabolic
subgroups of $G$ and let $\CP$ and $\CQ$ be the corresponding
generalized flag varieties. Let $u \in G^F_\uni$. Then we have
\begin{itemize}
\item[(i)]
$|\CP_u^F| = |\CQ_u^F|$;
\item[(ii)]
$|(\CP_u^0)^F| = |(\CQ_u^0)^F|$.
\end{itemize}
\end{cor}

\begin{rem}
\label{rem:Green-relations}
Assume as in Proposition \ref{prop1}.
In the special case  $P = G$ of Proposition \ref{prop1}(i),
we recover the following well-known identity
\[
\sum_{w \in W}  Q_{T_w}^G(u) = |W|,
\]
for any $u \in G_\uni^F$, e.g.\ see \cite[Prop.\ 5.8]{springer:green}.

Let $P$ be an $F$-stable parabolic subgroup of $G$
with $U = R_\uu(P)$ and let $\CP$ be the corresponding generalized
flag variety. Let $u \in G^F_\uni$.
Suppose that $G \cdot u \cap U = \varnothing$.
Then $|(\CP_u^0)^F| = 0$ and it follows from
Proposition \ref{prop1}(ii) that
\begin{equation}
\label{eq:zero1}
\sum_{w \in W_L} (-1)^{l(w)} Q_{T_w}^G(u) = 0.
\end{equation}
The special case $P = G$ (and $u \ne 1$)
of \eqref{eq:zero1} gives
\begin{equation}
\label{eq:zero2}
\sum_{w \in W} (-1)^{l(w)} Q_{T_w}^G(u) = 0,
\end{equation}
which is the formula in \cite[Prop.\ 5.2(ii)]{shoji}
for the sign character of $W$ and  $1 \ne u \in G_\uni^F$
(under the Springer correspondence, the sign representation of $W$
corresponds to the trivial unipotent class in $G$).
Conversely, we can recover
\eqref{eq:zero1} from the version of \eqref{eq:zero2} for ``$G = L$''
by applying the Harish-Chandra induction functor $R_L^G$.

Further, the special
case $P = G$ and $u = 1$ in Proposition \ref{prop1}(ii) gives
\begin{equation*}
\label{eq:zero3}
\frac{1}{|W|}\sum_{w \in W} (-1)^{l(w)} Q_{T_w}^G(1) = q^{\dim \CB},
\end{equation*}
which is just the formula in \cite[Prop.\ 5.2(i)]{shoji} again for the
sign character of $W$.
\end{rem}

\begin{rem}
\label{rem:uniform}
Note that $\CP_x$ and $\CP_x^0$ are defined for any $x \in G$
(though $\CP_x^0$ is empty unless $x$ is unipotent and
$G\cdot x \cap R_\uu(P) \ne \varnothing$).
It follows from the proofs of Lemmas  \ref{lem3}--\ref{lem5} that the functions
$G^F \to \ZZ_{\ge 0}$ given by $x \mapsto |\CP_x^F|$
and $x \mapsto |(\CP_x^0)^F|$ are \emph{uniform}, that is,
they are linear combinations of generalized Deligne--Lusztig
characters $R_{T_w}^G(\theta)$.
\end{rem}

\subsection{Polynomial properties for $|\CP_u^F|$ and $|(\CP_u^0)^F|$}
\label{s:poly} In this subsection we use the axiomatic setup given
in \S \ref{s:setup}. We fix the data $\Delta = (\Psi,\Pi,F_0,J)$,
where $\Psi = (X, \Phi, \check X, \check\Phi)$. Then for each prime
power $q$ we have a connected reductive algebraic group $G$ defined over
$\FF_q$ with Frobenius morphism $F$, and we have a parabolic
subgroup $P = P_J$ of $G$ which is $F$-stable. We assume that
$X/\ZZ\Phi$ is torsion free and that $q$ is a power of a good prime
for $\Phi$. We use the convention discussed in Remark \ref{R:varyq}
to vary $q$.

We recall that, with the above assumptions, the results of \S
\ref{s:classes} say that the parameterization of the unipotent
$G^F$-conjugacy classes does not depend on $q$. More specifically,
there is a finite set $E$ and finite groups $A_e$ ($e \in E$), such
that the unipotent $G^F$-conjugacy classes are parameterized by
pairs $(e,c)$ where $e \in E$ and $c$ is a conjugacy class of $A_e$.
This allows us to use the convention in Remark \ref{R:varyq2} to
take a representative $u = u_{e,c}$ of the conjugacy class with
label $(e,c)$ that possibly changes as we vary $q$.

We can now state the main result of this section.

\begin{thm}
\label{thmmain}
Fix the data $\Delta = (\Psi,\Pi,F_0,J)$, where
$\Psi = (X, \Phi, \check X, \check\Phi)$.  Assume that $X/\ZZ\Phi$
is torsion free and let $q$ be a power of a good prime for $\Phi$.
Let $G$, $F$ and $P$ be the connected reductive group, Frobenius morphism and
$F$-stable parabolic subgroup of $G$ determined by $\Delta$ and $q$.
Fix a label $(e,c)$ of a unipotent $G^F$-conjugacy class
and let $u \in G_\uni^F$ be a representative of this class.
\begin{itemize}
\item[(i)]
Suppose that $G$ does not admit a simple component of type $E_8$.
Then there exist $g_u(z), h_u(z)
\in \ZZ[z]$ such that $|\CP_u^F| = g_u(q)$ and $|(\CP_u^0)^F| =
h_u(q)$.
\item[(ii)]
Suppose that $G$ admits a simple component of type $E_8$.
Then there exist $g_u^i(z), h_u^i(z)
\in \ZZ[z]$ ($i = \pm 1$) such that $|\CP_u^F| = g_u^i(q)$ and
$|(\CP_u^0)^F| = h_u^i(q)$ when $q$ is congruent to $i$ modulo $3$.
\end{itemize}
\end{thm}

\begin{proof}
First suppose that $G$ does not have a simple component of type $E_8$.
In this case, the discussion in \S \ref{s:dlgenchars}
tells us that the values of the Green functions $Q_{T_w}^G(u)$ are given by
polynomials in $q$.  Now it follows from Proposition
\ref{prop1}(i) that there exists $g_u(z) \in
\QQ[z]$ such that $|\CP_u^F| = g_u(q)$.  That this
polynomial has integer coefficients follows from Proposition
\ref{p:intcoeffs} in the next subsection.

Observe that, as $L^F$ is a finite reductive group, the factor
$|L^F|_p$ in the formula for $|(\CP_u^0)^F|$
in Proposition \ref{prop1}(ii) is the size
of a Sylow $p$-subgroup of $L^F$.  So by \cite[Prop.\
3.19]{dignemichel}, we have $|L^F|_p = q^N$, where $N$ is the number
of positive roots in the root system of $L$.  Now it follows from
Proposition \ref{prop1}(ii) that there exists $k_u(z) \in \QQ[z]$
such that $|(\CP_u^0)^F| = k_u(q)/q^N$.  Since this holds for all
applicable $q$, Lemma \ref{L:rattopoly} tells us that $h_u(z) =
k_u(z)/z^N$ is a polynomial in $z$. Again, the fact that $h_u(z)$
has integer coefficients follows from Proposition \ref{p:intcoeffs}
in the next subsection.

In case $G$ has a simple component of type $E_8$, the proposition is
proved in the same way.  In this case we use Remark \ref{R:E8polys} to ensure
we have integer coefficients.
\end{proof}

\begin{rem}
\label{r:rem15}
We note that in Theorem \ref{thmmain}(ii)
the need to consider polynomials depending on the congruence of $q$ modulo
$3$ only arises in certain special cases, see Remark \ref{r:funnyE8}.
\end{rem}

\begin{rem}
\label{r:r16} Remove the assumption that $X/\ZZ\Phi$ is torsion
free, so that the centre of $G$ is possibly disconnected. Let $\hat
G$ be the semisimple adjoint group of the same type as $G$, and let
$\phi : G \to \hat G$ be the natural map as in \S \ref{s:classes}.
Then it is straightforward to see that, for $u \in G_\uni^F$, we
have $|\CP_u^F| = |\hat\CP_{\phi(u)}^F|$ and $|(\CP_u^0)^F| =
|(\hat\CP_{\phi(u)}^0)^F|$.  Therefore, with a certain formulation,
a version of Theorem \ref{thmmain} holds when $G$ has disconnected
centre.
\end{rem}

Using the {\tt chevie} package in GAP3 (\cite{gap}) along with some
code provided by M.~Geck, we can calculate the polynomials $g_u(z)$
and $h_u(z)$ from Theorem \ref{thmmain}. We illustrate this with some
examples below.

\begin{exmps} \label{ex:polys}
(i). Suppose $G$ is of type $D_5$. Let $u \in G_\uni^F$ lie in the minimal
unipotent $G$-class.  We note that the centralizer of $u$ in $G$ is
connected, so the elements of $G_\uni^F$ in the minimal unipotent
$G$-class form a single $G^F$-class.  We have
\begin{align*}
|\CB_u^F| & = 5q^{13} + 24q^{12} + 60q^{11} + 106q^{10} + 145q^9 +
161q^8 + 150q^7 + 120q^ 6 + 85q^5 + 54q^4 \\ & + 30q^3 + 14q^2 + 5q
+ 1
\end{align*}
and for $P$ a minimal parabolic subgroup we have
\begin{align*}
|\CP_u^F| & = q^{13} + 8q^{12} + 24q^{11} + 47q^{10} + 70q^9 + 83q^8
+ 82q^7 + 69q^6 + 51q^ 5 + 34q^4 \\ & + 20q^3 + 10q^2 + 4q + 1,
\\
|(\CP_u^0)^F| & = 4q^{12} + 16q^{11} + 36q^{10} + 59q^9 + 75q^8 +
78q^7 + 68q^6 + 51q^5 + 34q^4 + 20q^3 \\ & + 10q^2 + 4q + 1 .
\end{align*}

(ii). Now let $G$ be of type $G_2$ and
consider the unipotent $G$-class with Bala--Carter label
$G_2(a_1)$.  This is the subregular unipotent class in $G$.
Let $u_1 \in G_\uni^F$ be a split element in this
class. Then the component group of $C_G(u_1)$ is the symmetric group
$S_3$ on $3$ letters.  Letting $u_2 \in (G \cdot u_1)^F$ be in the
$G^F$-class corresponding to the conjugacy class of transpositions
in $S_3$ and $u_3 \in (G \cdot u_1)^F$ be in the $G^F$-class
corresponding to the conjugacy class of 3-cycles.  Then we have
\[
|\CB_{u_1}^F|  = 4q + 1, \qquad
|\CB_{u_2}^F|  = 2q + 1, \text{ and } \qquad
|\CB_{u_3}^F|  = q + 1.
\]

(iii). Finally we consider $G$ of type $F_4$.  Let $P$ be the standard
parabolic subgroup of $G$ with Levi factor of type $A_2$ corresponding to the
long roots of $\Pi$.  For $u$ the split representative in the
conjugacy class with Bala--Carter label $F_4(a_3)$, we have
$$
|(\CP_u^0)^F| = 16q + 4.
$$
For $u$ split in the conjugacy class with Bala--Carter label
$C_3(a_1)$, we have
$$
|(\CP_u^0)^F| = 4q^2 + 10q + 3.
$$
We note that in these cases the constant term in the polynomial is
not 1; in Remark \ref{rem:betti}(ii), we give an explanation for this.
\end{exmps}

\begin{exmp}
\label{ex:subreg} Let $G$ be simple and split over $\FF_q$, and let
$u \in G^F_\uni$ be a split subregular element of $G$. Let $r$ be
the rank of $G$. In \cite[III 3.10 Thm.\ 2]{steinberg} R.~Steinberg
showed that $\CB_u$ is a union of $r$ projective lines with an
intersection pattern determined by the Dynkin diagram of $G$; in
this case $\CB_u$ is called a \emph{Dynkin curve}. It follows that
$|\CB_u^F| = rq+1$, when $G$ is simply laced, $|\CB_u^F| =
(2r-1)q+1$, $rq+1$ in case $G$ is of type $B_r$, $C_r$ respectively,
and $|\CB_u^F| = (r+2)q+1$ for $G$ of type $G_2$ and $F_4$.
\end{exmp}

\subsection{Geometric interpretation of the polynomials}
\label{s:betti} In this subsection we show that the polynomials
$g_u(z)$ and $h_u(z)$ from Theorem \ref{thmmain} have integer
coefficients. Further, in case $G$ and $u$ are split we are able to
interpret the coefficients as Betti numbers of the relevant variety.
For simplicity, we only consider the case where $G$ does not have a
simple component of type $E_8$; analogous results for the case where
$G$ is of type $E_8$ can be proved similarly, see Remark
\ref{R:E8polys} for precise statements. We use the axiomatic setup
from \S \ref{s:setup}, so we fix the data $\Delta =
(\Psi,\Pi,F_0,J)$, where $\Psi = (X,\Phi,\check X, \check \Phi)$.
Then for any prime power $q$ we have corresponding $G$, $F$ and $P$.
Throughout this subsection we only consider powers of good primes
for $\Phi$. Further, we assume that $X/\ZZ\Phi$ is torsion free, so
that the centre of $G$ is connected; we note that Remark \ref{r:r16}
means that this latter assumption gives no essential loss in
generality.

We begin by considering the case where $F_0$ is the identity; this
corresponds to $G$ being split over $\FF_q$.  Further, we begin by only
considering unipotent $G^F$-conjugacy classes that contain a split
element $u$ (see \S \ref{s:classes}), and we concentrate on the
varieties $\CP_u$ and polynomials $g_u(z)$.

With the above conditions we can consider
$\CP_u$ as a variety defined over $\ZZ$,
so that for each prime power $q$ the definition over $\FF_q$ is obtained by
changing base. In particular, if we fix a prime $p$, then $\CP_u$ can
be viewed as a variety defined over $\FF_p$ and, for each power $q = p^s$
of $p$, the fixed points $\CP_u^F$ of the corresponding
Frobenius morphism are precisely the $\FF_q$-rational points of $\CP_u$ with
respect to the $\FF_p$-structure.

The discussion in the paragraph above,
Lemma \ref{L:intcoeffs} and the
purity of the variety $\CP_u$ (as discussed in \S \ref{s:purity})
imply the first part of Proposition \ref{p:betti} below; the second
part holds for the same reasons.

\begin{prop}
\label{p:betti}
Fix the data $\Delta = (\Psi,\Pi,F_0,J)$, where
$\Psi = (X,\Phi,\check X,\check \Phi)$. Suppose $\Phi$ does not have
an irreducible component of type $E_8$, the automorphism $F_0$ is
the identity and $X/\ZZ\Phi$ is torsion free.  Let $q$ be a power of
a good prime for $\Phi$.  Let $G$, $F$ and $P$ be the reductive
group, Frobenius morphism and $F$-stable parabolic subgroup of $G$
determined by $\Delta$ and $q$. Fix a split element $u \in G_\uni^F$.
\begin{enumerate}
\item[(i)]
Let $g_u(z)$ be as in Theorem \ref{thmmain}. Then
\[
g_u(z) = \sum_{i=0}^{d_u} \dim H^{2i}(\CP_u,\barQl) z^i,
\]
where $d_u = \dim \CP_u$.
\item[(ii)]
Let $h_u(z)$ be as in Theorem \ref{thmmain}. Then
\[
h_u(z) = \sum_{i=0}^{d_u^0} \dim H^{2i}(\CP^0_u,\barQl) z^i,
\]
where $d_u^0 = \dim \CP^0_u$.
\end{enumerate}
Further, for $i$ odd we have $H^i(\CP_u,\barQl) = \{0\} =
H^i(\CP_u^0,\barQl)$ and all the eigenvalues of $F$ on
$H^{2i}(\CP_u,\barQl)$ and $H^{2i}(\CP_u^0,\barQl)$ are $q^i$ for
all $i$. In particular, the Euler characteristics of $\CP_u$ and
$\CP^0_u$ are given by $g_u(1)$ and $h_u(1)$, respectively.
\end{prop}

The special case $P = B$ in Proposition \ref{p:betti}
is well-known, e.g.\ see \cite[Cor.\ 6.4]{shoji}.

\begin{rems}
\label{rem:betti}
(i). We note that the values $\dim H^i(\CP_u,\barQl)$
and $\dim H^i(\CP^0_u,\barQl)$ are the \emph{Betti numbers} of the
underlying varieties $\CP_u$ and $\CP_u^0$, respectively. Thus,
Proposition \ref{prop1} gives an effective way to calculate these
invariants.  For example, consider $G$ of type $D_5$ and $u$ in the
minimal unipotent $G$-conjugacy class and $P$ a minimal parabolic
subgroup. The expressions given in Example \ref{ex:polys}(i) tell us,
for instance, that $\dim H^{10}(\CP_u,\barQl) = 51$ and
$\dim H^{20}(\CP_u^0,\barQl) = 36$.

(ii).
It follows from \cite[Cor.\ 2.4]{steinberg2} that $\CP_u$ is
connected. Thus the constant term of $g_u(z)$ is $1$. In contrast,
it is not always the case that $\CP_u^0$ is connected, this follows
from the values of $h_u(z)$ for $G$ of type $F_4$ given in Example
\ref{ex:polys}(iii); another example is given in
\cite[II 11.8(b)]{spaltensteinbook}.
\end{rems}

\begin{rem}
\label{rem:betti-assoc} Let $P$ and $Q$ be associated $F$-stable
parabolic subgroups of $G$ and let $\CP$ and $\CQ$ be the
corresponding generalized flag varieties. Then it follows from
Corollary \ref{cor6} and Proposition \ref{p:betti} that for $u \in
G_\uni$, the Betti numbers of $\CP_u$ and $\CQ_u$
(respectively $\CP_u^0$ and $\CQ_u^0$) coincide. We note that this
has been observed by Borho--MacPherson (\cite[Cor.\ 3.7]{borhomac}).
In particular, $\dim H^{2d_u}(\CP_u, \barQl) = \dim H^{2d_u}(\CQ_u,
\barQl)$ implies that $\CP_u$ and $\CQ_u$ have the same number of
irreducible components of maximal dimension. We note that it is not
always the case that $\CP_u$ and $\CQ_u$ have the same number of
irreducible components in lower dimensions; see \cite[II
11.5]{spaltensteinbook} for an example in case $G = \GL_4$.

Analogous statements hold if we replace $\CP_u$ and
$\CQ_u$ by $\CP_u^0$ and $\CQ_u^0$; see \cite[II 11.6]{spaltensteinbook}
for an explicit example where the number of irreducible components
differs. We note however that for $G = \GL_n$, the varieties $\CP_u^0$
are equidimensional, see the final corollary in \cite{spaltenstein0};
also see \cite[II Prop.\ 5.15]{spaltensteinbook}.

Note that in contrast to $\CP_u$ and $\CP_u^0$, the fixed point variety
$\CB_u$ is always equidimensional, \cite[II Prop.\ 1.12]{spaltensteinbook}.
\end{rem}

\begin{rem}
\label{rem:paving} Next we discuss some special cases where there is
an alternative way to see that $|\CP_u^F|$ is a polynomial in $q$
with coefficients independent of $q$. Let $G$ and $u \in G_\uni^F$
be split. Suppose that $\CP_u$ admits a stratification into a finite
union of locally closed subvarieties, so that the decomposition does
not depend on $q$ and each stratum is isomorphic to some affine
space and is defined over ${\mathbb F}_q$; such a decomposition is
called a {\em paving by affine spaces}. It then follows that
$|\CP_u^F|$ is a polynomial in $q$. Further, the coefficient of
$q^i$ in this polynomial is the number of strata isomorphic to
affine $i$-space $\AA^i$; by Proposition \ref{p:betti} this
coefficient is also the $(2i)$th Betti number of $\CP_u$. Below we
give some examples when such a stratification is known to exist.

(i).   In case $u = 1$, the Bruhat decomposition of $G$ implies that
there is a paving of $\CB_1 = \CB$ by affine spaces, and one deduces
that $\dim H^{2i}(\CB,\barQl)$ is the number of elements of $W$ of
length $i$. Likewise, the Bruhat decomposition can be used to obtain a
stratification of $\CP$ and
to interpret the $(2i)$th Betti number of $\CP$ as the number of elements
of length $i$ in the set of distinguished coset representatives in $W/W_L$.

(ii) Suppose $G$ is simple. It is known in all cases apart from
$E_7$ and $E_8$ that $\CB_u$ admits a paving by affine spaces for
any $u \in G_\uni$, see \cite[Thm.\ 3.9]{deconcinietal}, \cite[II
Prop.\ 5.9]{spaltensteinbook}, and \cite[Thm.\ 2]{spaltenstein1}.
See also \cite[Ch.\ 11]{jantzen}. We note that the treatment in
\cite{deconcinietal} only applies in characteristic $0$. In case $u$
is in the subregular class of $G_\uni$, $\CB_u$ is a Dynkin curve
(as mentioned in Example \ref{ex:subreg}). Therefore, $\CB_u$ is the
disjoint union of affine lines and points with no restriction on the
type of $G$. Also, if $u$ is a root element in $G_\uni$, then
$\CB_u$ admits an affine paving, irrespective of the type of $G$,
see \cite[Prop.\ 7.11]{springer:green}.

(iii). Now suppose $G = \GL_n$ and $P$ is a parabolic
subgroup of $G$ with associated generalized flag variety $\CP$.  In
\cite{shimomura}, N.~Shimomura showed that $\CP_u$ admits a paving
by affine spaces for any $P$ and any $u$ in $G_\uni$.
\end{rem}

In view of the special situations discussed in Remark \ref{rem:paving}
it is natural to ask whether in general $\CP_u$ and $\CP_u^0$
admit pavings by affine spaces.
In particular, this is consistent with the vanishing of the odd
$l$-adic cohomology of $\CP_u$ and $\CP_u^0$,
as shown in Proposition \ref{p:betti}.

\smallskip

Now we consider the situation when $u$ is not split (still assuming that $G$
is split).  So we fix a label $(e,c)$ of a unipotent $G^F$-conjugacy
class and let $u$ be a representative of this conjugacy class;
we use Remark \ref{R:varyq2} to vary $q$.  We recall that $u$ is split
if the permutation action of $F$ on the irreducible components of $\CB_u$
is trivial.

Since there are finitely many irreducible components of $\CB_u$,
there is a minimal positive integer $s$ such that the action of
$F^s$ on the set of irreducible components of $\CB_u$ is trivial.
(We note that $s$ is the order of an element of the conjugacy class
$c$ in the group $A_e$, but as we do not require this fact we do not
give a proof of it.)  Therefore, considered as an element of
$G^{F^s}$ we have that $u$ is split.  Hence, by Proposition
\ref{p:betti} above, the $l$-adic cohomology groups vanish in odd
degrees, and the eigenvalues of $F^s$ on the $l$-adic cohomology
groups $H^{2i}(\CP_u,\barQl)$ are all $q^{is}$ for $i = 0,\dots,d_u$.
Therefore, the eigenvalues of $F$ on $H^{2i}(\CP_u,\barQl)$ are of
the form $\zeta^j q^i$, where $\zeta$ is a primitive $s$th root of
unity and $j = 0,\dots,s-1$.  Let $n_{i,j}$ be the number of times
$\zeta^jq^i$ ($j = 0,\dots,s-1$) occurs as an eigenvalue of $F$ on
$H^i(\CP_u,\barQl)$.  Using the Grothendieck trace formula we get
\begin{equation}
\label{eq:puf}
|\CP_u^F| = \sum_{i=0}^{d_u} \sum_{j=0}^{s-1} n_{i,j} \zeta^j q^i.
\end{equation}
Since \eqref{eq:puf} holds for all $q$ and $|\CP_u^F| =
g_u(q)$, we see that the $n_{i,j}$ do not depend on $q$.  It follows
that the coefficients of $g_u(z)$ are algebraic integers.  Since these
coefficients are also rational numbers, we must have $g_u(z) \in \ZZ[z]$.

Next we consider the case where we no longer assume that $F_0$ is
the identity, i.e.\ that $G$ is not split over $\FF_q$.  The
automorphism $F_0$ is of finite order, say $r$.  Then $F^r$ is a
Frobenius morphism on $G$ corresponding to a split
$\FF_{q^r}$-structure.  The discussion in the previous paragraph
then tells us that the eigenvalues of $F^r$ on
$H^{2i}(\CP_u,\barQl)$ are of the form $\zeta^j q^{ir}$, where
$\zeta$ is some root of unity.  It follows that the eigenvalues of
$F$ on $H^{2i}(\CP_u,\barQl)$ are of the form $\xi^j q^i$, where
$\xi$ is some root of unity. Now an analogue of \eqref{eq:puf} holds
and we can deduce as before that $g_u(z) \in \ZZ[z]$.

Proposition \ref{p:intcoeffs} completes the proof of Theorem
\ref{thmmain}. We have proved the statement about $g_u(z)$ and the
assertion about $h_u(z)$ is proved in the same way.

\begin{prop}
\label{p:intcoeffs} Fix the data $\Delta = (\Psi,\Pi,F_0,J)$, where
$\Psi = (X, \Phi, \check X, \check\Phi)$.  Assume that $\Phi$ does
not have a simple component of type $E_8$ and $X/\ZZ\Phi$ is torsion
free.  Let $q$ be a power of a good prime for $\Phi$. Let $G$, $F$
and $P$ be the connected reductive group, Frobenius morphism and $F$-stable
parabolic subgroup of $G$ determined by $\Delta$ and $q$.
Fix a label $(e,c)$ of a
unipotent $G^F$-conjugacy class and let $u \in G_\uni^F$ be a
representative of this class.

Let $g_u(z), h_u(z)$ be as in Theorem \ref{thmmain}.  Then we have
$g_u(z),h_u(z) \in \ZZ[z]$.
\end{prop}

\begin{rem} \label{R:E8polys}
Suppose $G$ is split and has a simple component of type $E_8$.  Then
for $u \in G_\uni^F$ split, the formulas in Proposition
\ref{p:betti} hold for $g_u^1(z)$ and $h_u^1(z)$. Further, the
analogue of Proposition \ref{p:intcoeffs} holds without the
assumption that $G$ and $u$ are split, i.e.\ that $g_u^i(z),h_u^i(z)
\in \ZZ[z]$ ($i = \pm 1$). As mentioned in Remark \ref{r:rem15}, we
only actually need to consider two polynomials depending on the
congruence of $q$ modulo 3 in certain special cases.
\end{rem}

\section{Unipotent conjugacy in finite groups of Lie type}
\label{s:alperin}

In this section we apply the results of the previous sections in
order to prove our generalization (Theorem \ref{thm:main}) of the
main theorem in Alperin's note \cite{alperin}.

\subsection{A counting argument for finite groups}
\label{s:notation} For the first part of this subsection let $G$ be
a finite group and let $H$ be a subgroup of $G$.  In Lemma
\ref{lem1} below we give a formula for the number of $H$-conjugacy
classes in $G$.  This lemma was proved for a particular case in
\cite{alperin}, but the proof applies generally without change; we
include the proof for the reader's convenience. Before stating the
lemma we need to introduce some notation.

For $x \in G$ we write
\begin{equation*}
f_H^G(x) = |\{{}^gH \mid x \in {}^gH, g \in G\}|
\end{equation*}
for the number of conjugates of $H$ in $G$ containing $x$.  By
$k(H,G)$ we denote the number of $H$-conjugacy classes in $G$ and we
let $\CR = \CR(G,H)$ be a set of representatives of the distinct
$G$-conjugacy classes that intersect $H$.

\begin{lem}
\label{lem1}
The number of $H$-conjugacy classes in $G$ is given by
\begin{equation*}
k(H,G) = \frac{|N_G(H)|}{|H|}\sum_{x \in \CR}f_H^G(x).
\end{equation*}
\end{lem}

\begin{proof}
Counting the set of pairs $\{(y,{}^gH) \mid y \in G \cdot x, g \in G, y
\in {}^gH\}$ in two different ways, we obtain $|G:N_G(H)|\cdot |G \cdot
x \cap H| = |G \cdot  x|\cdot  f_H^G(x)$. It follows that
\begin{equation*}
|G \cdot  x \cap H| = \frac{|G \cdot  x|}{|G:N_G(H)|} \cdot  f_H^G(x).
\end{equation*}
Now Burnside's counting formula gives
\begin{equation*}
k(H,G) = \frac{1}{|H|} \sum_{x \in H} |C_G(x)| = \frac{1}{|H|}
\sum_{x \in \CR} |G \cdot  x \cap H|\cdot |C_G(x)|.
\end{equation*}
Thus it follows that
\[
k(H,G) = \frac{1}{|H|} \sum_{x \in \CR} |G \cdot  x| \cdot
\frac{|N_G(H)|}{|G|}\cdot \frac{|G|}{|G \cdot x|}\cdot f_H^G(x) =
\frac{|N_G(H)|}{|H|}\sum_{x \in \CR}f_H^G(x),
\]
as claimed.
\end{proof}

\begin{rem} \label{R:allreps}
If the conjugacy class of $x \in G$ does not intersect $H$, then we
have $f_H^G(x) = 0$.  Therefore, it does no harm in the formula in
Lemma \ref{lem1} to sum over representatives of all conjugacy
classes of $G$, i.e.\ to take $\CR = \CR(G)$ to be a set of
representatives of all the conjugacy classes of $G$.
\end{rem}

Now suppose $G$ is a connected reductive group defined over $\FF_q$
with corresponding Frobenius morphism $F$ of $G$.  Let $H$ be a closed
$F$-stable subgroup of $G$.  In this situation we write $f_H^G$
rather than $f_{H^F}^{G^F}$ to shorten notation; we note that this
notation has already been used in \S \ref{Sfixedpoint}.  We can
apply Lemma \ref{lem1} in the case ``$G = G^F$'' and ``$H = H^F$''
and we obtain the formula
\[
k(H^F,G^F) = \frac{|N_G(H^F)|}{|H^F|}\sum_{x \in \CR}f_H^G(x),
\]
where $\CR = \CR(G^F,H^F)$ is a set of representatives of all
$G^F$-conjugacy classes that intersect $H^F$; or using Remark
\ref{R:allreps} we may take $\CR = \CR(G^F)$ to be a set of
representatives of all $G^F$-conjugacy classes.

Lastly in this subsection, we remark on a connection between the
functions $f_H^G$ and so called fixed point ratios for finite group
actions; for the importance of the latter see for instance
\cite{LLS} and the references therein.

\begin{rem}
\label{r:fpr} Let $G$ be a finite group and let $H$ be a subgroup.
Consider the transitive action of $G$ on the cosets $G/H$. For $x
\in G$, the \emph{fixed point ratio} of $x$ in this action is given
by
\[
\fpr(x, G/H) = \frac{\fix_{G/H}(x)}{|G/H|} = \frac{|G\cdot x \cap
H|}{|G \cdot x|},
\]
where $\fix_{G/H}(x)$ denotes the number of fixed points of $x$ on
$G/H$, see \cite{LLS}. In particular, by the counting argument in
the proof of Lemma \ref{lem1} we get
\[
\fpr(x, G/H) = \frac{|N_G(H)|}{|G|} f_H^G(x).
\]
\end{rem}

\subsection{Unipotent conjugacy in finite groups of Lie type}
For the start of this section we do not require the axiomatic setup
of \S \ref{s:setup}, but we do use the notation of \S \ref{s:notn}.
So $G$ is a connected reductive algebraic group defined over $\FF_q$ with
corresponding Frobenius morphism $F$.  The parabolic subgroup $P$ of $G$ is
$F$-stable with unipotent radical $U$.  We have a maximally split
maximal torus $T$ of $G$ and $L$ is the unique Levi subgroup of $P$
containing $T$.

The next proposition is proved using Lemmas \ref{lem5} and \ref{lem1}.

\begin{prop}
\label{propUG}
The number of $U^F$-conjugacy classes in $G^F$ is
given by
\[
k(U^F,G^F)  = |L^F|_{p'} \cdot \frac{1}{|W_L|}\sum_{w \in W_L}
(-1)^{l(w)} \sum_{u \in \CR}Q_{T_w}^G(u),
\]
where $\CR = \CR(G^F,U^F)$ is a set of representatives of the
conjugacy classes of $G^F$ that intersect $U^F$.
\end{prop}

\begin{proof}
We have $N_{G^F}(U^F) = P^F$ and therefore $|N_{G^F}(U^F)|/|U^F| =
|L^F|$. Thus it follows from Lemmas \ref{lem5} and \ref{lem1} that
\begin{align}
\label{eq:kug}
k(U^F,G^F)
& = |L^F| \sum_{u \in \CR}f_U^G(u) \notag \\
& = |L^F| \sum_{u \in \CR} \frac{1}{|L^F|_{p}} \cdot
\frac{1}{|W_L|}\sum_{w \in W_L}
(-1)^{l(w)} Q_{T_w}^G(u) \\
& =  |L^F|_{p'} \cdot \frac{1}{|W_L|}\sum_{w \in W_L} (-1)^{l(w)}
\sum_{u \in \CR}Q_{T_w}^G(u), \notag
\end{align}
as claimed.
\end{proof}

For the remainder of this subsection we use the axiomatic setup of
\S \ref{s:setup}. We fix the data $\Delta = (\Psi,\Pi,F_0,J)$, where
$\Psi = (X,\Phi,\check X,\check \Phi)$.  Then for each prime power
$q$ we have a connected reductive algebraic group $G$ defined over $\FF_q$
with Frobenius morphism $F$ and we have the $F$-stable parabolic subgroup $P =
P_J$ of $G$. We denote the unipotent radical of $P$ by $U$ and we
write $L$ for the unique Levi subgroup of $P$ containing $T$. In
order to state the results in this subsection we use the convention
discussed in Remark \ref{R:varyq} to vary $q$.

We now state and prove the principal result of this section; it
generalizes the main theorem of \cite{alperin}, which is the special
case $G = \GL_n$ and $P = B$ of Theorem \ref{thm:main}.

\begin{thm}
\label{thm:main} Fix the data $\Delta = (\Psi,\Pi,F_0,J)$, where
$\Psi = (X, \Phi, \check X, \check\Phi)$. For a prime power $q$, let
$G$, $F$ and $P$ be the connected reductive group, Frobenius morphism and
$F$-stable parabolic subgroup of $G$ determined by $\Delta$ and
let $U = R_\uu(P)$.
Assume that $X/\ZZ\Phi$ is torsion free and that $q$ is a power of a
good prime for $\Phi$.
\begin{itemize}
\item[(i)]
Suppose that $G$ does not have a simple component of type $E_8$.
Then there exists $m(z) \in \ZZ[z]$ such that $k(U^F,G^F) = m(q)$.
\item[(ii)]
Suppose that $G$ has a simple component of type $E_8$. Then there exist
$m^i(z) \in \ZZ[z]$ ($i = \pm 1$), such that $k(U^F,G^F) = m^i(q)$,
when $q$ is congruent to $i$ modulo $3$.
\end{itemize}
\end{thm}

\begin{proof}
The first equality in \eqref{eq:kug}
says that
\[
k(U^F,G^F) = |L^F| \sum_{u \in \CR}f_U^G(u).
\]
Since $L^F$ is a finite reductive group, the factor $|L^F|$ is a
polynomial in $q$, e.g., see \cite[p.\ 75]{carter}. By Lemma
\ref{lem3}(ii) and Theorem \ref{thmmain} each of the summands
$f_{U}^G(u)$ is a polynomial in $q$. With the assumption that
$X/\ZZ\Phi$ is torsion free, it follows from Proposition
\ref{prop:classes2} that the set $\CR = \CR(G^F,U^F)$ is independent
of $q$, where we use the convention of Remark \ref{R:varyq2}.
\end{proof}

\begin{rem}
\label{rem:e8maxpar}
By Remark \ref{r:funnyE8}, if $G$ has a simple
component of type $E_8$ and $P$ is an $F$-stable parabolic subgroup of
$G$, then the polynomial $k(U^F,G^F)$ in $q$ in Theorem
\ref{thm:main}(ii) only depends on the residue class of $q$ modulo
$3$ in case $U$ meets the unipotent class with Bala--Carter label
$E_8(b_6)$. For $G$ simple of type $E_8$, it is straightforward to
calculate the parabolic subgroups of $G$ for which this occurs; for
example this class does not meet the unipotent radical of any of the
maximal parabolic subgroups of $G$.
\end{rem}

It was already remarked that the formula in Lemma \ref{lem5} for
$f_U^G(u)$ does not depend on the choice of an $F$-stable parabolic
subgroup $P$ of $G$ containing $L$ as a Levi subgroup.  This enabled
us to deduce Corollary \ref{cor6}.  Similarly, we can observe the
following corollary.  We can argue as in Remark \ref{R:allreps} to
sum over all unipotent $G^F$-conjugacy classes in the formula for
$k(U^F,G^F)$ in Proposition \ref{propUG}; alternatively one can note
that associated parabolic subgroups have the same Richardson class,
so that the unipotent classes of $G$ meeting their unipotent
radicals coincide.

\begin{cor}
\label{cor2}
Assume as in Theorem \ref{thm:main}.
Let $P$ and $Q$ be associated $F$-stable parabolic
subgroups of $G$ with unipotent radicals $U$ and $V$ respectively.
Then
\[
k(U^F,G^F) = k(V^F,G^F).
\]
\end{cor}

\begin{rem}
\label{rem:comvar} Assume as in Theorem \ref{thm:main}. The
\emph{commuting variety of $U$ and $G$} is the closed subvariety of
$U \times G$ defined by
\[
\CC(U,G) = \{(u,g) \in U \times G \mid ug = gu \}.
\]
Since $U$ is $F$-stable, $F$ acts on $\CC(U,G)$ and we can consider the
$F$-fixed points and we have $\CC(U,G)^F = \CC(U^F,G^F)$.
The Burnside counting formula tells us that
\[
|\CC(U^F,G^F)| = q^{\dim U} k(U^F,G^F) = q^{\dim U} m(q),
\]
where $m(z) \in \ZZ[z]$ is as in Theorem \ref{thm:main}. In
particular, $|\CC(U^F,G^F)|$ is given by a polynomial in $q$ with
integer coefficients independent of $q$. (If $G$ admits a simple
component of type $E_8$, there are two polynomials depending on the
congruence class of $q$ modulo $3$.)  We note that the coefficients
of the polynomial $m(z)$ can be negative (see Examples
\ref{ex:polys2}). This means that the variety $\CC(U,G)$ is not, in
general, pure and therefore we cannot interpret the coefficients as
Betti numbers as in Lemma \ref{L:intcoeffs}.
%
\end{rem}

Using the {\tt chevie} package in GAP3 (\cite{gap}) along with some
code provided by M.~Geck, we are able to explicitly calculate the
polynomials $m(z)$ in Theorem \ref{thm:main}. We illustrate this with some
examples below.

\begin{exmps}
\label{ex:polys2}
(i). In Table \ref{Tab:GL} below we give the polynomials for $k(U^F,G^F)$
in case $G = \PGL_n$ and $P = B$ for $n = 2,\dots,10$.  In this case
we can take $U = \U_n$ to be the group of upper unitriangular
matrices.  We note that it follows from the first equality in
\eqref{eq:kug} that $k(\U_n(q),\GL_n(q)) =
(q-1)k(\U_n(q),\PGL_n(q))$, because $f_{\U_n}^{\GL_n}(u) =
f_{\U_n}^{\PGL_n}(u)$ for any $u \in \U_n$.

\begin{table}[h!tb]
\renewcommand{\arraystretch}{1.5}
\begin{tabular}{|l|p{400pt}|}
\hline $n$ & $k(\U_n(q),\PGL_n(q))$ \\
\hline\hline 2 & $q^2 + q - 2$ \\
\hline
3 & $q^5 + q^3 - 3q^2 - 2q + 3$ \\
\hline
4 & $q^9 - q^7 + 2q^6 - 3q^5 - 2q^4 - 4q^3 + 9q^2 + 3q - 5$ \\
\hline 5 & $q^{14} - q^{12} - q^{11} + 3q^{10} - 5q^9 + 2q^8 - 9q^6
+ 11q^5
+ 10q^4 - q^3 - 13q^2 - 4q + 7$ \\
\hline 6 & $q^{20} - q^{18} - q^{17} - q^{16} + 5q^{15} - 6q^{14} -
3q^{13} + 12q^{12} + 3q^{11} - 40q^{10} + 44q^9 - 7q^8 + 5q^7 + 3q^6
- 5q^5 - 38q^4 +
11q^3 + 24q^2 + 5q - 11$ \\
\hline 7 & $q^{27} - q^{25} - q^{24} - q^{23} + 6q^{21} - 6q^{20} -
4q^{19} + 3q^{18} + 16q^{17} + 2q^{16} - 45q^{15} + 9q^{14} +
65q^{13} - 36q^{12} - 47q^{11} + 118q^{10} - 130q^9 + 80q^8 - 85q^7
+ 25q^6 + 34q^5 + 46q^4 -
27q^3 - 31q^2 - 6q + 15$ \\
\hline 8 & $q^{35} - q^{33} - q^{32} - q^{31} + 8q^{28} - 7q^{27} -
5q^{26} + 3q^{25} + 2q^{24} + 27q^{23} - 8q^{22} - 76q^{21} +
66q^{20} + 9q^{19} + 8q^{18} - 96q^{17} + 109q^{16} + 56q^{15} -
73q^{14} - 266q^{13} + 357q^{12} + 93q^{11} - 530q^{10} + 278q^9 +
253q^8 - 153q^7 - 52q^6 + 11q^5 -
96q^4 + 51q^3 + 48q^2 + 7q - 22$ \\
\hline 9 & $q^{44} - q^{42} - q^{41} - q^{40} + q^{37} + 9q^{36} -
8q^{35} - 6q^{34} + 3q^{33} + q^{32} + 8q^{31} + 25q^{30} - 2q^{29}
- 113q^{28} + 49q^{27} + 107q^{26} - 60q^{25} + 81q^{24} - 326q^{23}
+ 97q^{22} + 702q^{21} - 603q^{20} - 446q^{19} + 337q^{18} + 760q^
{17} - 869q^{16} + 491q^{15} - 957q^{14} + 1063q^{13} - 142q^{12} +
123q^{11} - 939q^{10} + 1130q^9 - 622q^8 - 255q^7 + 429q^6 - 60q^5 +
124q^4 - 92q^3 -
60q^2 - 8q + 30$ \\
\hline  10 & $q^{54} - q^{52} - q^{51} - q^{50} + q^{47} + q^{46} +
10q^{45} - 9q^{44} - 7q^{43} + 2q^{42} + q^{41} + 9q^{40} - 3q^{39}
+ 41q^{38} - 12q^{37} - 144q^{36} + 61q^{35} + 77q^{34} + 89q^{33} -
90q^{32} - 5q^{31} - 189q^{30} - 109q^{29} + 561q^{28} + 256q^{27} -
746q^{26} - 50q^{25} - 1070q^{24} + 3249q^{23} - 682q^{22} -
4884q^{21} + 3467q^{20} + 5522q^{19} - 8703q^{18} + 757q^{17} +
5424q^{16} - 1423q^{15} - 1450q^{14} - 4812q^{13} + 10000q^{12} -
6872q^{11} + 726q^{10} + 1638q^9 - 555q^8 + 509q^ 7 - 858q^6 +
307q^5 - 222q^4 + 137q^3 + 85q^2 + 9q -
42$ \\
\hline
\end{tabular}
\medskip
\caption{$k(\U_n(q),\PGL_n(q))$} \label{Tab:GL}
\end{table}

(ii). Now suppose that $G$ is split of type $E_6$. For $P = B$, we have
\begin{align*}
m(q) & = q^{42} - q^{40} - q^{37} - q^{36} + q^{35} - q^{33} +
q^{32} + 8q^{31} - 2q^{30} - 6q^{29} - 2q^{28} - 6q^{27} + 18q^{26}
\\ & - 24q^{25} - 6q^{24} + 34q^{23} + 90q^{21} - 122q^{20} - 59q^{19}
- 69q^{18} + 101q^{17} + 304q^{16} - 265q^{15} \\ & + 234q^{14} +
48q^{13} - 1891q^{12} + 2658q^{11} - 220q^{10} - 1794q^9 + 818q^8 +
527q^7 - 521q^6 \\ & + 213q^5 + 11q^4 - 78q^3 - 17q^2 - 6q + 25,
\end{align*}
and if $P$ is a minimal parabolic
subgroup of $G$, then $k(U^F,G^F)$ is given by
\begin{align*}
m(q) & = q^{43} - q^{41} - q^{38} - q^{37} + q^{36} - q^{34} +
q^{33} + 7q^{32} - q^{31} - 5q^{30} - 2q^{29} - 6q^{28} + 12q^{27}
\\ & - 13q^{26} - 9q^{25} + 23q^{24} + 7q^{23} + 65q^{22} - 61q^{21}
-81q^{20} - 53q^{19} + 47q^{18} + 231q^{17} \\
& - 78q^{16} + 42q^{15}  + 79q^{14} - 987q^{13} + 748q^{12} +
931q^{11} - 1191q^{10} - 129q^9 + 563q^8 \\ & - 39q^7 - 117q^6 +
74q^5 - 37q^4 - 40q^3 - 3q^2 + 24q.
\end{align*}

(iii). Now consider the case $G$ is of type $E_6$ with twisted definition
over $\FF_q$, i.e.\ $G^F = {}^2E_6(q)$.  Then for a Borel subgroup
of $G$ we have
\begin{align*}
m(q) & = q^{42} - 8q^{41} + 31q^{40} - 80q^{39} + 160q^{38} -
273q^{37} + 423q^{36} - 615q^{35} + 848q^{34} - 1113q^{33} \\ & +
1401q^{32} - 1710q^{31} + 2039q^{30} - 2376q^{29} + 2699q^{28} -
2992q^{27} + 3258q^{26} - 3521q^{25} \\ & + 3821q^{24} - 4201q^{23}
+ 4698q^{22} - 5347q^{21} + 6155q^{20} - 7035q^{19} + 7831q^{18}  -
8471q^{17} \\ & + 9005q^{16} - 9454q^{15}  + 9764q^{14} - 9947q^{13}
+ 10114q^{12} - 10268q^{11} + 10141q^{10} - 9392q^9 \\ & + 7989q^8 -
6252q^7 + 4526q^6 - 2998q^5 + 1770q^ 4 - 893q^3 + 354q^2 - 94q + 12.
\end{align*}
For a minimal parabolic subgroup of $G$, we have
\begin{align*}
m(q) & = q^{43} - 6q^{42} + 17q^{41} - 32q^{40} + 48q^{39} -
65q^{38} + 85q^{37} - 107q^{36} + 126q^{35} - 139q^{34} \\ & +
149q^{33} - 160q^{32}
 + 169q^{31} - 168q^{30} + 155q^{29} - 138q^{28} + 127q^{27} - 130q^{26} +
154q^{25} \\
& - 200q^{24} + 264q^{23}  - 341q^{22} + 409q^{21} - 420q^{20} +
357q^{19} - 280q^{18} + 233q^{17} - 171q^{16} \\
& + 71q^{15}  - 4q^{14} + 18q^{13} - 29q^{12} - 107q^{11} +
410q^{10} - 718q^9 + 863q^8 - 834q^7 + 706q^6 \\
& - 531q^5 + 346q^4 - 182q^3 + 65q^2 - 11q.
\end{align*}
\end{exmps}

In case $G$ has disconnected centre, the parameterization of the
unipotent $G^F$-conjugacy classes depends on $q$, as discussed at
the end of \S \ref{s:classes}.  We also mentioned that this
dependence is only up to congruences on $q$. From this discussion,
Theorem \ref{thm:main} and Remark \ref{r:r16} it follows that
$k(U^F,G^F)$ is given by polynomials up to congruences on $q$. We
illustrate this in the following example.

\begin{exmp}
\label{ex:disc}
Let $G = \SL_3$. In this case, we have three $F$-stable unipotent
$G$-conjugacy classes, as mentioned in Example \ref{ex:SL3}. The
class corresponding to the partition $(1,1,1)$ is $\{1\}$, which
forms a single $G^F$-conjugacy class; one can also show that the
$F$-stable class corresponding to the partition $(2,1)$ is a single
$G^F$-class. As discussed in Example \ref{ex:SL3} the regular class
splits into three $G^F$-classes if $q \equiv 1 \mod 3$ and forms a
single $G^F$-class if $q \equiv 0,-1 \mod 3$. We consider the case
when $P = B$ is a Borel subgroup of $G$.
One can calculate the value $f_U^G(u)$ to be $1$, $2q+1$,
$(q^2 + q + 1)(q+1)$ for $u$ in the class labeled $(3)$, $(2,1)$,
$(1,1,1)$, respectively.  Using the formula given in the first
equation in \eqref{eq:kug}, we obtain
\begin{align*}
k(U^F,G^F) & = (q-1)^2(3 + (2q+1) + (q^2 + q + 1)(q+1)) \\
& = q^5 + q^3 - q^2 - 6q + 5,
\end{align*}
for $q \equiv 1 \mod 3$; and
\begin{align*}
k(U^F,G^F) & = (q-1)^2(1 + (2q+1) + (q^2 + q + 1)(q+1)) \\
& =  q^5 + q^3 - 3q^2 - 2q + 3,
\end{align*}
for $q \equiv 0,-1 \mod 3$.  So
$k(U^F,G^F)$ is given by two polynomials depending on congruences of
$q$ modulo $3$.
\end{exmp}

The formula in Proposition \ref{propUG} and Remark
\ref{rem:badchargreen} suggest that the values of $k(U^F,G^F)$
differ for good and bad primes.  Below we give an example where this
does occur.

\begin{exmp}
\label{ex:b2} Consider $G$ of type $C_2$ and $P = B$ a Borel
subgroup of $G$. If $q$ is a power of an odd prime, then we have
\[
k(U^F,G^F) = q^6 + 4q^3 - 8q^2 - 2q + 5
\]
while if $q$ is a power of $2$, then we get
\[
k(U^F,G^F) = q^6 + q^4 + 2q^3 - 6q^2 - 4q + 6.
\]
We are grateful to M.~Geck for the computation of these polynomials.
\end{exmp}

\begin{rem}
\label{rem:q-1} Suppose that $G$ is split and that $P$ be an
$F$-stable parabolic subgroup of $G$ such that every element in
$\CR$ is split. This situation is uncommon, but for example it is
the case if $G = \GL_n$.  Let $L$ be an $F$-stable Levi subgroup of $P$
and let $u \in \CR$. It follows from Lemma \ref{lem3}, Theorem
\ref{thmmain}, and Proposition \ref{p:betti} that $f_U^G(u)$ is a
polynomial in $q$ with non-negative integer coefficients. Since $L$
is reductive, it follows from the Bruhat decomposition of $L^F$ that
$|L^F|$ is a polynomial in $q-1$ with non-negative integer
coefficients, cf.~\cite[Prop.\ 3.19(ii)]{dignemichel}. By the first
equality in \eqref{eq:kug} we have $k(U^F,G^F) = |L^F| \sum_{u \in
\CR}f_U^G(u)$. Consequently, $k(U^F,G^F)$ is a  polynomial in $q-1$
with non-negative integer coefficients.

We observe that for any of the polynomials $k(U^F,G^F)$ that we have
calculated in this paper, we have $k(U^F,G^F) \in \NN[q-1]$. It
seems likely that this positivity behaviour always holds and it
would be interesting to know if there is a geometric explanation for
it.
\end{rem}

\begin{rem}
\label{rem:p}
Let $P$ be an $F$-stable parabolic subgroup of $G$.
Using Lemma \ref{lem1} and the fact that $P = N_G(P)$ we obtain
\begin{equation}
\label{eq:kPG}
k(P^F,G^F)  =  \sum_{x \in \CR} f_P^G(x),
\end{equation}
where $\CR = \CR(P^F,G^F)$ is a set of representatives of the
conjugacy classes of $G^F$ that intersect $P^F$.  From the proof of
Lemma \ref{lem4} we get an expression of $f_P^G(x)$ as a sum of
Deligne--Lusztig generalized characters, but this is a linear
combination of Green functions only in case $x$ is unipotent.
Further, the set $\CR$ over which the sum is taken depends heavily
on $q$, so it is non-trivial to determine whether $k(P^F,G^F)$ is
given by a polynomial in $q$.

In \cite{alperin} it is shown, using an expression as in
\eqref{eq:kPG}, that $k(B^F,G^F)$ is given by a polynomial in $q$,
in case $G = \GL_n$.  The proof of this depends on partitioning the
set $\CR(B^F,G^F)$ into a finite union $\CR(B^F,G^F) = \CR_1 \cup
\dots \cup \CR_t$ independent of $q$ such that: $f_B^G(x) =
f_B^G(y)$ if $x,y \in \CR_i$; and $|\CR_i|$ is given by a polynomial
in $q$. An inductive counting argument is used to show that
$f_B^G(x)$ is given by a polynomial in $q$, for $x \in G$.

We hope to use similar arguments to consider $k(P^F,G^F)$ in general.
In \cite{goodwinroehrle:parabolic}, we have proved that
$k(P^F,G^F)$ is a polynomial in $q$ in case $G = \GL_n$.  This is
achieved using a decomposition of the conjugacy classes of $G^F$
analogous the partition of $\CR(B^F,G^F)$ mentioned above.  Further,
we reduce the calculation of $f_P^G(x)$ to that of $f_Q^H(u)$ for
certain pseudo Levi subgroups $H = C_G(s)^\circ$ of $G$ and
parabolic subgroups $Q$ of $H$, where $x = su$ is the Jordan
decomposition of $x$.  We hope that similar arguments will apply to
arbitrary reductive $G$.  As the centre of a pseudo Levi subgroup of
$G$ need not be connected even if the centre of $G$ is connected,
one can only really hope to prove that $f_P^G(x)$ is given by
polynomials in $q$ up to congruences on $q$.



Finally, we observe that we readily obtain an analogue of
Corollary \ref{cor2}. Let $P$ and $Q$ be two associated $F$-stable
parabolic subgroups of $G$. Then we have
\begin{equation}
k(P^F,G^F) = k(Q^F,G^F).
\end{equation}
We argue as follows. We may assume that $P$ and $Q$ share a common
$F$-stable Levi subgroup $L$ say. Then the desired equality follows
from \eqref{eq:kPG},
and the equality $f_P^G(x) = R_L^G(1_L)(x) = f_Q^G(x)$ for any $x
\in G^F$, as $R_L^G(1_L)$ does not depend on the choice of an
$F$-stable parabolic subgroup admitting $L$ as a Levi subgroup.
\end{rem}

\appendix

\section{The $\GL_n$ case}

In this appendix, we give an elementary proof of the existence of
the polynomial $h_u(z)$ in Theorem \ref{thmmain} for the case $G =
\GL_n$ with standard definition over $\FF_q$. In this case existence
of $h_u(z)$ follows from the theory of Jordan normal forms, so one
obtains an elementary proof of Theorem \ref{thm:main}.  The proof is
purely combinatorial, beginning with two counting lemmas.

Let $V$ be an $n$-dimensional vector space over $\FF_q$ and let $W$
be an $m$-dimensional subspace of $V$.  Let $l \in \ZZ_{\ge 0}$ and
define
\[
\CS(V,W,l) = \{U \subseteq V \mid \dim U = l, \ U \cap W = \{0\}\}.
\]
We note that the order of $\CS(V,W,l)$ depends only on $n$, $m$ and $l$
(and not on the choice of $V$ and $W$).  As we let $q$ vary, the
order of $\CS(V,W,l)$ defines a function of $q$, which we denote by
$s(n,m,l)(q)$.

The following lemma is proved using an easy counting argument.

\begin{lem} \label{L:subspace}
We have
\[
s(n,m,l)(q) = \prod_{i=0}^{l-1} \frac{q^n - (q^m + q^i -
1)}{q^i-q^{i-1}}.
\]
In particular, $s(n,m,l)(q)$ is a polynomial in $q$.
\end{lem}

\begin{proof}
Suppose we have chosen subspace $U'$ of $V$
with $\dim U' = l-1$ and $U' \cap W = \{0\}$.
Then the number of ways of extending $U'$ to $U \in
\CS(V,W,l)$ is $(q^n - (q^m + q^{l-1} - 1))/(q^{l-1}-q^{l-2})$.  The
formula for $s(n,m,l)(q)$ now follows by induction.  We have that
$s(n,m,l)(q)$ is a polynomial in $q$ by Lemma \ref{L:rattopoly}
\end{proof}

We recall that a flag $F$ in $V$ is a sequence of subspaces $ \{0\} =
V_0 \sse \dots \sse V_r = V$; in this appendix we only consider proper flags
so $V_i \ne V_{i+1}$ for all $i$.  The dimension vector of $F$ is
the $r$-tuple $d = (d_1,\dots,d_r) \in \ZZ_{\ge 1}^r$, where $d_i = \dim V_i$.
Fix a flag $F$ in $V$
with dimension vector $d$ and let $e = (e_1,\dots,e_r) \in \ZZ_{\ge
0}^r$ be an $r$-tuple with $e_i \le d_i$ and $e_i \le e_{i+1}$
for each $i$. We define $\CT(F,e)$ to be
the set of subspaces $U$ of $V$ such that $\dim(U \cap V_i) = e_i$
for each $i$.  It is clear that the order of $\CT(F,e)$ depends only
on the vectors $d$ and $e$ and as $q$ varies, it defines a function $t(d,e)(q)$
of $q$. We prove the following lemma by induction using Lemma
\ref{L:subspace}.

\begin{lem} \label{L:flag}
$t(d,e)(q)$ is a polynomial in $q$.
\end{lem}

\begin{proof}
We work by induction on $r$.  The case $r = 0$ is trivial.
Therefore, we may assume that the number of ways of choosing $U' \sse
V_{r-1}$ such that $\dim(U' \cap V_i) = e_i$ for $i = 1,\dots,r-1$
is given by a polynomial in $q$.  Given such $U'$, the number of
ways of extending $U'$ to $U$ such that $U \cap V_{r-1} = U'$ and
$\dim(U \cap V_r) = e_r$ is the number of subspaces of
$V_r/U'$ of dimension $e_r - e_{r-1}$ that intersect $V_{r-1}/U'$
trivially.  The induction is now completed using Lemma
\ref{L:subspace}.
\end{proof}

\bigskip

We retain the notation that $V$ is an $n$-dimensional vector space
over $\FF_q$ and we define $G = \GL(V)$.  For a fixed dimension
vector $d = (d_1,\dots,d_r)$ ($d_r = n$, $d_{i-1} < d_i$), we define
$\Flag(V,d)$ to be the set of all flags of $V$ with dimension vector
$d$.

In this appendix we denote a partition $\pi$ of $n$ by the $n$-tuple
$(\pi_1,\dots,\pi_n)$ of non-negative integers such that $\sum_{i =
1}^n i\pi_i = n$. We denote by $\Pi(n)$ the set of all partitions of
$n$. We recall that the unipotent conjugacy classes in $G$ are
parameterized by the partitions of $n$ (by the theory of Jordan
normal forms). Given a partition $\pi$ of $n$, let $x_\pi \in G$ be
a representative of the conjugacy class corresponding to $\pi$; we
write $y_\pi = x_\pi - 1$.  For any unipotent element $x \in G$ we
write $\Par(x)$ for the element of $\Pi(n)$ to which the conjugacy
class of $x$ corresponds.

For a partition $\pi$ and dimension vector $d$ we define
\[
\CF(\pi,d) = \{F \in \Flag(V,d) \mid y_\pi(V_i) \sse V_{i-1} \text{
for all $i$}\}.
\]
As $q$ varies the order of $\CF(\pi,d)$ defines a function of $q$,
which we denote by $f(\pi,d)(q)$; note that this function depends only
on $\pi$ and not the choice of $x_\pi$.

Given a flag $F \in \Flag(V,d)$ we define $P_F$ to be the stabilizer
of $F$ in $G$, i.e.\
\[
P_F = \{x \in G \mid x(V_i) = V_i \text{ for all $i$}\}.
\]
It is well-known that $P_F$ is a parabolic subgroup of $G$ and that all
parabolic subgroups of $G$ occur in this way (for some $d$ and $F$).
The unipotent radical $U_F$ of $P_F$ is given by
\[
U_F = \{x \in G \mid (x-1)(V_i) \sse V_{i-1} \text{ for all $i$}\}.
\]
It is clear that for any $F \in \Flag(V,d)$ and any $\pi \in \Pi(n)$
there is a bijection between $\CF(\pi,d)$ and
\[
\{{}^gU_F \mid x_\pi \in {}^gU_F,\ g \in G\}.
\]
Therefore, the existence of the polynomial
$h_u(z)$ in Theorem \ref{thmmain} is equivalent to proving that
$f(\pi,d)(q)$ is a polynomial in $q$.

\smallskip

If $x \in G$ and $W$ is a subspace of $V$ such that $x(W) = W$, then
we write $x^W$ for the element of $\GL(V/W)$ induced by $x$. Let $c
\in \ZZ_{\ge 1}$ and let $\pi'$ be a partition of $n - c$. Define
\[
\CW(\pi,\pi') = \{W \sse V \mid \dim W = c,\ y_\pi(W) = \{0\},\
\Par(x_\pi^W) = \pi' \}.
\]
As $q$ varies, the order of $\CW(\pi,\pi')$ defines a function
$w(\pi,\pi')(q)$ of $q$, which only depends on $\pi$ and $\pi'$
(and not on the choice of $x_\pi$).

It is clear that
\[
f(\pi,d)(q) = \sum_{\pi' \in \Pi(n-d_1)}
w(\pi,\pi')(q)f(\pi',d')(q),
\]
where $d' = (d_2,\dots,d_r)$.  So by induction we just have to show
that $w(\pi,\pi')(q)$ is a polynomial in $q$, for all $\pi,\pi'$.

Define $K_\pi = \ker y_\pi$ and for each $j \in \ZZ_{\ge 0}$ define
$J_{\pi,j} = K_\pi \cap  y_\pi^j(V)$. Clearly, a subspace of $V$ is
killed by $y_\pi$ if and only if it is contained in $K_\pi$.  Let $v
\in K_\pi$ and let $j \in \ZZ_{\ge 0}$ be such that $v \in J_{\pi,j}
\setminus J_{\pi,j-1}$. It is easy to check that the endomorphism
that $x_\pi$ induces on $V/\langle v \rangle$ corresponds to the
partition $\pi^*$, where $\pi^*_j = \pi_j - 1$, $\pi^*_{j-1} =
\pi_{j-1} + 1$ and $\pi^*_i = \pi_i$ for all other $i$.

Let $\pi' \in \Pi(n-d_1)$ and for $i = 1,\dots,n$ we define $b_i =
\sum_{j=i}^n \pi_j - \pi'_j$. The observations in the previous
paragraph imply that $W \in \CW(\pi,\pi')$ if and only if $W \sse
K_{\pi}$ and $\dim(W \cap J_{\pi,i}) = b_i$ for each $i$. It now
follows from Lemma \ref{L:flag} that $w(\pi,\pi')(q)$ is a
polynomial in $q$.

This completes the proof of existence of the polynomial $h_u(z)$ in
Theorem \ref{thmmain} and therefore the proof of Theorem
\ref{thm:main} for $G = \GL_n$.

We have not justified that $h_u(z) \in \ZZ[z]$, this reduces to
checking that the polynomial $s(n,m,l)(q)$ from Lemma
\ref{L:subspace} has integer coefficients; we leave the details to
the reader.

Further, we note that it is far from evident from this proof
that for associated parabolic subgroups of $G = \GL_n$
we get the same polynomials,
cf.\ Corollaries \ref{cor6} and \ref{cor2}.


\bigskip

{\bf Acknowledgments}: We are grateful to J.~M.~Douglass and M.~Geck
for helpful comments. Part of the research for this paper was
carried out while both authors were staying at the Mathematical
Research Institute Oberwolfach supported by the ``Research in
Pairs'' programme at Oberwolfach. Also, part of the paper was
written during a visit of both authors to the Max-Planck-Institute
for Mathematics in Bonn. We thank both institutions for their
support. In addition, part of this research was funded by an EPSRC
grant.  The first author thanks New College, Oxford for financial
support and hospitality whilst the research was completed.



\bigskip

\begin{thebibliography}{88}

\bibitem{alperin}
J.~L.~Alperin,
\emph{Unipotent conjugacy in general linear groups},
Comm. Algebra \textbf{34} (2006), no.\ 3, 889--891.

\bibitem{beyspalt}
W.~M.~Beynon and N.~Spaltenstein,
\emph{Green functions of finite Chevalley groups of type
$E\sb{n}$ $(n=6,\,7,\,8)$},
J.\ Algebra \textbf{88} (1984), no.\ 2, 584--614.

\bibitem{borhomac}
W.~Borho and R.~MacPherson,
\emph{Partial resolutions of nilpotent
varieties}, Analysis and topology on singular spaces, II, III
(Luminy, 1981), 23--74, Ast\'erisque, 101-102, Soc.\ Math.\ France,
Paris, 1983.

\bibitem{carter}
R.~W.~Carter,
\emph{Finite groups of Lie type. Conjugacy classes and
complex characters}, Pure and Applied Mathematics,
New York, 1985.

\bibitem{deconcinietal}
C.~De Concini, G.~Lusztig, and C.~Procesi,
\emph{Homology of the zero-set of a nilpotent vector field on a flag manifold},
J. Amer. Math. Soc. \textbf{1} (1988), no. 1, 15--34.

\bibitem{dignemichel}
F.~Digne and J.~Michel,
\emph{Representations of finite groups of Lie type},
London Mathematical Society Student Texts {\bf 21}, Cambridge
University Press, Cambridge, 1991.

\bibitem{gap}
The GAP group,
\emph{GAP -- Groups, Algorithms, and Programming -- version 3 release
4 patchlevel 4}, Lehrstuhl D f\"ur Mathematik, Rheinisch
Westf\"alische Technische Hochschule, Aachen, Germany, 1997.

\bibitem{geck}
M.~Geck,
\emph{On the average values of the irreducible characters
of finite groups of Lie type on geometric unipotent classes}, Doc.\
Math.\ \textbf{1} (1996), No.\ 15, 293--317.

\bibitem{goodwinroehrle:parabolic}
S.~M.~Goodwin and G.~R\"ohrle,
\emph{Parabolic conjugacy in general linear groups},
J. Algebraic Combin., \textbf{27}, no.~1, (2008), 99--111.

\bibitem{higman}
G.~Higman,
\emph{Enumerating $p$-groups. I. Inequalities}, Proc.\
London Math.\ Soc.\ (3) \textbf{10} (1960) 24--30.

\bibitem{jantzen}
J.~C.~Jantzen, \emph{Nilpotent orbits in representation theory}, Lie
theory, 1--211, Progr. Math., 228, Birkh\"auser Boston, Boston, MA,
2004.

\bibitem{johnrich}
D.~S.~Johnston and R.~W.~Richardson, \emph{Conjugacy classes in
parabolic subgroups of semisimple algebraic groups. II}, Bull.\
London Math.\ Soc.\ \textbf{9} (1977), no.\ 3, 245--250.

\bibitem{kawanaka}
N.~Kawanaka,
\emph{Generalized Gelfand-Graev representations of
exceptional simple algebraic groups over a finite field. I},
Invent.\ Math.\ \textbf{84} (1986), no.\ 3, 575--616.

\bibitem{LLS}
R.~Lawther, M.W.~Liebeck, and G.M.~Seitz,
\emph{Fixed point ratios in actions of finite exceptional groups of Lie type},
Pacific J. Math. \textbf{205} (2002), no. 2, 393--464.

\bibitem{LusztigV}
G.~Lusztig,
\emph{Character sheaves. V}.
Adv. in Math. \textbf{61} (1986), no. 2, 103--155.

\bibitem{mcnsom}
G.~J.~McNinch and E.~Sommers,
\emph{Component groups of unipotent centralizers in good characteristic},
J.\ Algebra {\bf 260} (2003), no.\ 1, 323--337.

\bibitem{mizuno}
K.~Mizuno,
\emph{The conjugate classes of unipotent elements of the Chevalley groups
$E\sb{7}$ and $E\sb{8}$}. Tokyo J. Math. \textbf{3} (1980), no. 2, 391--461.

\bibitem{pomerening}
K.~Pommerening,
\emph{\"Uber die unipotenten Klassen reduktiver Gruppen}, J.\
Algebra {\bf 49} (1977), no.\ 2, 525--536.

\bibitem{premet}
A.~Premet, \emph{Nilpotent orbits in good characteristic and the
Kempf-Rousseau theory}, J.\ Algebra {\bf 260} (2003), no.\ 1,
338--366.

\bibitem{robinson}
G.~R.~Robinson, \emph{Counting conjugacy classes of unitriangular
groups associated to finite-dimensional algebras}, J.\ Group Theory
\textbf{1} (1998), no.\ 3, 271--274.

\bibitem{shimomura}
N.~Shimomura,
\emph{A theorem on the fixed point set of a
unipotent transformation on the flag manifold},
J. Math. Soc. Japan \textbf{32} (1980), no. 1, 55--64.

\bibitem{shoji}
T.~Shoji, \emph{Green functions of reductive groups over a finite
field}, The Arcata Conference on Representations of Finite Groups
(Arcata, Calif., 1986), 289--301, Proc.\ Sympos.\ Pure Math., 47,
Part 1, Amer.\ Math.\ Soc., Providence, RI, 1987.

\bibitem{spaltenstein0}
N.~Spaltenstein,
\emph{The fixed point set of a unipotent transformation
on the flag manifold},
Nederl. Akad. Wetensch. Proc. Ser. A \textbf{79}
(1976), no. 5, 452--456.

\bibitem{spaltensteinbook}
\bysame, \emph{Classes unipotentes et sous-groupes de Borel},
Lecture Notes in Mathematics, \textbf{946} Springer-Verlag,
Berlin-New York, 1982.

\bibitem{spaltenstein1}
\bysame,
\emph{On unipotent and nilpotent elements of groups of type $E\sb{6}$},
J. London Math. Soc. (2) \textbf{27} (1983), no. 3, 413--420.

\bibitem{spaltenstein}
\bysame,
\emph{On the reflection representation in
Springer's theory}, Comment.\ Math.\ Helv.\  \textbf{66}  (1991),
no.\ 4, 618--636.

\bibitem{springer:green}
T.~A.~Springer,
\emph{Trigonometric sums, Green functions of finite groups and
representations of Weyl groups}  Invent. Math.  \textbf{36}  (1976), 173--207.

\bibitem{springer1}
\bysame,
\emph{A purity result for fixed point varieties in
flag manifolds}, J.\ Fac.\ Sci.\ Univ.\ Tokyo Sect. IA Math.
\textbf{31} (1984), no.\ 2, 271--282.


\bibitem{springersteinberg}
T.~A.~Springer and R.~Steinberg,
\emph{Conjugacy classes}.
Seminar on Algebraic Groups and Related Finite Groups
(The Institute for Advanced Study, Princeton, N.J., 1968/69)
pp. 167--266
in Lecture Notes in Mathematics, \textbf{131}
Springer-Verlag, Berlin-New York, 1970.



\bibitem{steinberg}
R.~Steinberg, \emph{Conjugacy classes in algebraic groups},
Lecture Notes in Mathematics, Vol.\ 366, Springer-Verlag, Berlin-New
York, 1974.

\bibitem{steinberg2}
\bysame,
\emph{On the desingularization of the unipotent variety},
Invent. Math. \textbf{36} (1976), 209--224.

\bibitem{thompson}
J.~Thompson, \emph{$k(\U_n(F_q))$}, Preprint,
{\tt{http://www.math.ufl.edu/fac/thompson.html}}.

\bibitem{veralopezarregi}
A.~Vera-L\'opez and J.~M.~Arregi, \emph{Conjugacy classes in unitriangular
matrices}, Linear Algebra Appl.\  \textbf{370}
(2003), 85--124.

\end{thebibliography}
\end{document}